\documentclass[a4paper,fleqn]{cas-sc}

\usepackage[authoryear,longnamesfirst]{natbib}

\usepackage{amsmath}
\usepackage{amssymb}
\usepackage{amsfonts}
\usepackage{amsthm}
\usepackage{mathtools}
\usepackage{mathrsfs}
\usepackage{booktabs}
\usepackage{multirow}
\usepackage{tikz}

\newtheoremstyle{assumptionstyle}%
  {}
  {}
  {\itshape}
  {}
  {\bfseries}
  {}
  { }
  {\thmname{#1}\thmnumber{ #2}\thmnote{ (#3)}}

\theoremstyle{assumptionstyle}

\theoremstyle{plain}
\newtheorem{thm}{Theorem}[section]
\newtheorem{lem}{Lemma}[section]

\newtheorem{rem}{Remark}

\renewcommand{\P}{\mathbb{P}}
\renewcommand{\Pr}{\mathrm{P}}

\newcommand{\E}{\mathbb{E}}

\newcommand{\cov}{\text{Cov}}

\newcommand{\perm}{\mathrm{perm}}
\newcommand{\xvec}{\mbox{\bf x}}

\newcommand{\Xvec}{\mbox{\bf X}}
\newcommand{\Yvec}{\mbox{\bf Y}}
\newcommand{\Zvec}{\mbox{\bf Z}}

\newcommand{\wvec}{\mbox{\bf w}}

\newcommand{\zvec}{\mbox{\bf z}}
\newcommand{\yvec}{\mbox{\bf y}}
\newcommand{\thetavec}{\mbox{\boldmath $\theta$}}

\newcommand{\sigmat}{\mbox{\boldmath $\Sigma$}}
\newcommand{\muvec}{\mbox{\boldmath $\mu$}}

\DeclareMathOperator{\Var}{Var}

\begin{document}
\let\WriteBookmarks\relax
\def\floatpagepagefraction{1}
\def\textpagefraction{.001}

\shorttitle{High-dimensional two-sample testing}
\shortauthors{Ghosh et al.}

\title[mode=title]{High-dimensional two-sample testing: From optimal transport to optimal derangements}

\author[1]{Abhradipta Ghosh}
\author[1]{Vaibhav Sherkar}
\author[2]{Bilol Banerjee}
\author[1]{Anil K. Ghosh}

\cormark[1]
\ead{anilkghosh@gmail.com} \cortext[1]{Corresponding Author}

\affiliation[1]{
  organization={Theoretical Statistics and Mathematics Unit, Indian Statistical Institute},
  city={Kolkata},
  country={India}
}

\affiliation[2]{
  organization={Department of Statistics and Data Science, National University of Singapore},
  city={Singapore},
  country={Singapore}
}

\begin{abstract}
In recent years, optimal transport has emerged as a useful framework for constructing distribution-free multivariate two-sample tests by transporting observations to a reference distribution under a suitable transportation cost. However, existing optimal transport-based tests can perform poorly in high dimension, low sample size (HDLSS) settings, and their performance depends critically on the choice of the reference distribution and transportation cost. We investigate this phenomenon and show that suitable choices 
can substantially improve their performance.



Motivated by this analysis, we introduce a class of two-sample tests based on optimal derangements that do not require an external reference distribution.
The minimum distance pairing test of \cite{MinDist} arises as a special case when the optimization is restricted to a subclass of derangements, thereby revealing a connection between this method and optimal transport. We propose the optimal derangement (ODER) test and its distribution-free variant, the ODER-SHC test, obtained by restricting the minimization to derangements consisting of a single Hamiltonian cycle.
We establish their consistency under both HDLSS and classical large-sample asymptotic regimes. We also develop corresponding procedures based on the \(\ell_1\) distance and establish their high-dimensional consistency under conditions involving coordinate-wise energy distances. These \(\ell_1\) distance based tests can also discriminate between two high-dimensional distributions with identical first two moments but different marginal distributions.  
Analyses of benchmark datasets and extensive simulation studies demonstrate the performance of the proposed methods across a variety of alternatives and data settings.

\end{abstract}



\maketitle

\section{Introduction}
In a two-sample problem, we test for the equality of two $d$-dimensional distributions $F$ and $G$ based on two sets of independent observations $\{{\bf x}_1,{\bf x}_2,\ldots,{\bf x}_m\}$ 
and $\{{\bf y}_1,{\bf y}_2,\ldots,{\bf y}_n\}$ from $F$ and $G$, respectively. This is a fundamental problem in statistics, and a wide variety of tests are available for it. 
A summary of classical nonparametric tests for univariate data can be obtained in  
\cite{hollander2013nonparametric} and \cite{gibbons2014nonparametric}.
Several tests are available for multivariate data as well. Classical tests for the two-sample location problem include those based on coordinate-wise signs and ranks \citep[see, e.g.,][]{puri1971nonparametric}, spatial signs and ranks \citep[see, e.g.,][]{oja2010multivariate} and interdirections
\citep[see, e.g.,][]{randles1989distribution,hallin2002optimal}. 
However, these tests are generally not applicable when the dimension of the data exceeds the pooled sample size. To address high-dimensional location problems, several test procedures have been proposed \citep[see, e.g.,][]{bai1996effect,chen2010two,aoshima2018two}. Some projection-based tests \citep[see, e.g.,][]{ghosh2016distribution,wei2016direction,srivastava2016raptt}
are also available for this purpose. 
\cite{li2012two}, \cite{cai2013two} and \cite{ding2024two} proposed two-sample tests for equality of covariance matrices of two high-dimensional distributions. 

Several tests have been developed for the general two-sample problem as well. \cite{BARINGHAUS2004,BARINGHAUS2010} proposed a two-sample test based on energy distances. 
Similar tests were also considered by \cite{szekely2013energy} and \cite{biswas2014nonparametric}. However, these distance-based tests are sensitive to outliers. 
\cite{Kernel} constructed a kernel two-sample test based on maximum mean discrepancy (MMD), but its performance depends critically on the choice of the kernel function and the associated bandwidth. Several graph-based tests have also been proposed. \cite{MST} developed a runs test based on minimum spanning tree (MST). \cite{Neighbors1986} and \cite{Neighbors1988} proposed tests based on nearest neighbors. However, these tests often yield poor performance for high-dimensional data, especially when the scale difference between the two distributions dominates their location difference \citep[see, e.g.,][]{biswas2014nonparametric}. \cite{Chen02012017} proposed some modified tests to take care of these issues. Other notable tests that can be used for high-dimensional two-sample problems include \cite{liu2011triangle,mondal2015high, sarkar2018some,pan2018ball,tsukada2019high, sarkar2020some1} and \cite{Wass2020}. 
However, unlike univariate nonparametric tests, none of the above-mentioned tests are distribution-free. 
\cite{Run2014} developed a distribution-free run test, but its implementation involves solving an NP-hard optimization problem \citep{garey1979computers}. Moreover, its large-sample behavior has yet to be investigated in the literature. 
\cite{MinDist} proposed a distribution-free two-sample test based on minimum distance pairing (MDP). However, this test often performs poorly in high dimensions \citep[see, e.g.,][]{sarkar2020some1}.

Recently, \cite{OTTest} proposed a general recipe based on optimal transport 
to construct distribution-free tests for multivariate data. 
In particular, they considered the rank energy test, a rank version of the energy
test \citep{BARINGHAUS2004},  where the ranks are obtained using optimal transport. For this purpose, they used the square of the Euclidean distance as the transportation cost and $U[0,1]^d$, the uniform distribution on the unit hypercube, as the reference distribution. This test has large-sample consistency under appropriate regularity conditions, and can also be used for high-dimensional data even when the dimension is larger than the sample size.
So, naturally, one would be curious to know how this test performs in high dimensions. 

Our theoretical and empirical investigation in the next section shows that while the rank energy test based on the conventional choice of the reference distribution and the transportation cost \citep[as in][]{OTTest} often leads to poor performance in high dimension, low sample size (HDLSS) situations, an intelligent choice can substantially increase its power. One can also 
construct a simpler distribution-free test based on transport allocation count with similar high-dimensional behavior. 
However, the performance of these tests depends critically on the choice of the reference distribution, for which a suitable specification may not always be available for the data at hand. 
We address this issue in Section 3, where we introduce a class of tests based on optimal derangements. These tests do not require any external reference distribution for their implementations. Interestingly, \cite{MinDist}'s test based on MDP can be derived as a special case when the optimization is restricted to a particular subclass of derangements. 
We propose two tests, one based on unrestricted optimization and the other by restricting the optimization to the subclass of derangements
consisting of a single Hamiltonian cycle. In the latter case, the resulting test turns out to be distribution-free.
We investigate the HDLSS behavior of the proposed tests under appropriate assumptions. 
We also develop corresponding tests
based on the $\ell_1$ distance and study their high-dimensional behavior. These tests can also
discriminate between two high-dimensional distributions differing outside the first two
moments. Some benchmark datasets are analyzed
in Section 4 to compare the empirical performance of the proposed tests with some state-of-the-art methods. Some simulated datasets are also analyzed in Section 5 for further evaluation of their performance. Section 6 investigates the large-sample behavior of the proposed tests. In particular, we derive their asymptotic null distribution and prove their large-sample consistency under mild regularity conditions. Finally, Section 7 contains a brief summary of the work and ends with a discussion on some possible directions for future research. All proofs and mathematical details are given in the Appendix.

\section {Tests based on optimal transport and their high-dimensional behavior}

The idea of optimal transport \citep[see, e.g.,][]{OT,
 hallin2022measure} dates back to \cite{monge1781memoire}. Given two probability measures \(\mu_0\) and \(\nu_0\), the Monge formulation seeks a transformation that maps \(\mu_0\) to \(\nu_0\) while minimizing a prescribed transportation cost. Since such a transformation need not exist, \cite{kantorovich1948problem} addressed this issue by considering probability measures on the product space with marginals \(\mu_0\) and \(\nu_0\), and minimizing the expected transportation cost over all such measures. The resulting minimizer is called an optimal transport plan. 

\cite{OTTest} used optimal transport to construct distribution-free two-sample tests for multivariate data.  
Let $\{{\bf x}_1,{\bf x}_2,\ldots,{\bf x}_m\}$ 
and $\{{\bf y}_1,{\bf y}_2,\ldots,{\bf y}_n\}$ be independent observations from $F$ and $G$, respectively, and let $H_{m+n}$ denote the empirical distribution of the pooled sample $\{\zvec_1=\xvec_1,\ldots,\zvec_m=\xvec_m,\zvec_{m+1}=\yvec_1,\ldots,\zvec_{m+n}=\yvec_n\}$. 
They transported \(H_{m+n}\) to an empirical version of the uniform distribution on $[0,1]^d$, using $C(\zvec,\wvec)=\|\zvec-\wvec\|^2$, the squared Euclidean distance, as the transportation cost. 
If $\{\wvec_1,\wvec_2,\ldots,\wvec_{m+n}\}$ are generated from $U[0,1]^d$, the optimal transport amounts to finding a permutation \(\pi\) minimizing
the total transportation cost 
${\cal C}(\pi)=\sum_{i=1}^{m+n}\|{\zvec}_i-{\wvec}_{\pi(i)}\|^2$.
The assigned point \(\mathbf{w}_{\pi(i)}\) ($i=1,2,\ldots,m+n$) serves as the multivariate rank of $\zvec_i$. Under the null hypothesis \({\cal H}_0:F=G\), these ranks have a distribution determined only by the uniform reference sample, yielding distribution-free tests. 
\cite{OTTest} proposed a rank version of the energy statistic based on these ranks, referred to as the rank energy test. 

Note that this test remains applicable in high dimensions, even when the dimension exceeds the pooled sample size. To illustrate the performance of the rank energy test in the HDLSS setting, consider the following two examples:

{\bf Example A:} Two normal distributions
$N_d({\bf 0}_d,{\bf I}_d)$ and 
$N_d(0.5 {\bf 1}_d,{\bf I}_d)$ differ only in their locations.

{\bf Example B:} Two normal distributions
$N_d({\bf 0}_d,{\bf I}_d)$ and
$N_d({\bf 0}_d,2{\bf I}_d)$
differ only in their scales.

Here ${\bf 0}_d$ and ${\bf 1}_d$ are $d$-dimensional vectors with all elements equal to $0$ and $1$, respectively, ${\bf I}_d$
is the $d \times d$ identity matrix and $N_d(\muvec,\sigmat)$ denotes the $d$-dimensional normal distribution with mean $\muvec$ and dispersion matrix $\sigmat$.

We carried out our experiment for $m=n=20$ and different values of $d$ ($d=2^i$ for $i=1,2,\ldots,11$). In each case, the experiment was repeated 500 times, and the empirical power of the test was computed as the proportion of times it rejected ${\cal H}_0$. Figure 1 shows that in both examples,
the rank energy test performed
poorly in high dimensions. In the scale problem, its power remained close to the nominal level of 0.05 for all values of $d$. In the location problem, it had somewhat better performance in small dimensions, but its power decreased rapidly as $d$ increased. 

\begin{figure}[htbp]
    \centering
    \includegraphics[height=1.7in,width=0.4\textwidth]{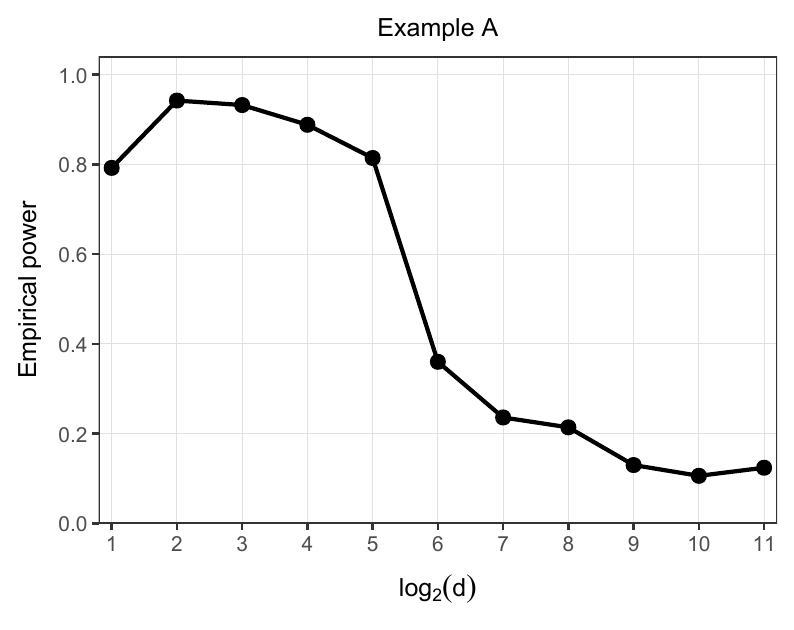}
    \includegraphics[height=1.7in,width=0.4\textwidth]{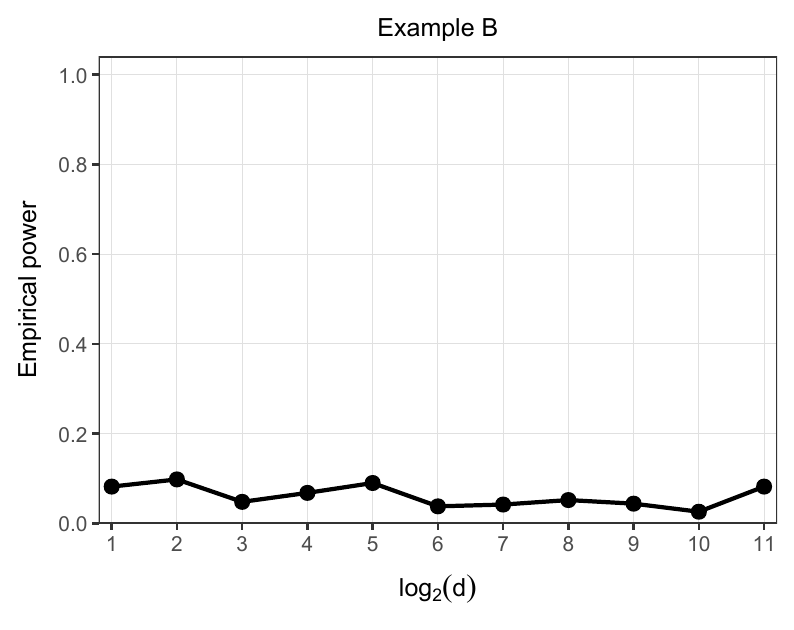}

    \vspace{-0.1in}
    \caption{Powers of the rank energy test in Examples A and B}
    \label{fig:rank-energy}
\end{figure}

This behavior can be explained by the concentration of transportation costs. 
If ${\bf X}\sim N_d({\bf 0}_d,{\bf I}_d)$, ${\bf Y} \sim N_d(\mu{\bf 1}_d,\sigma^2{\bf I}_d)$ and 
${\bf U} \sim U[0,1]^d$ are independent; using the strong law of large numbers (SLLN), one can show that ${d}^{-1}\|{\bf X}-{\bf U}\|^2 \stackrel{a.s.}{\rightarrow}  4/3$ and ${d}^{-1}\|{\bf Y}-{\bf U}\|^2 \stackrel{a.s.}{\rightarrow} (\mu-1/2)^2+\sigma^2+{1/12}$  
as $d \rightarrow \infty$. 
Consequently, ${\cal C}(\pi)/d$ has asymptotically the same limit for every permutation $\pi$. Thus, in high dimensions, the optimal permutation becomes essentially arbitrary, and the distribution of the resulting rank energy statistic turns out to be almost the same as its null permutation distribution. This explains the loss of power observed in Figure 1.



\begin{figure}[!b]
    \centering
    \includegraphics[height=1.7in,width=0.4\textwidth]{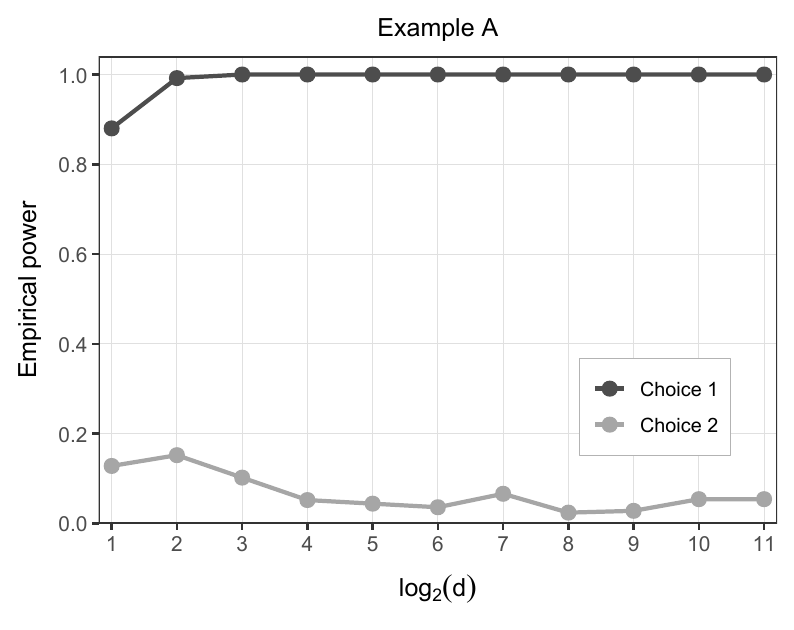}
    \includegraphics[height=1.7in,width=0.4\textwidth]{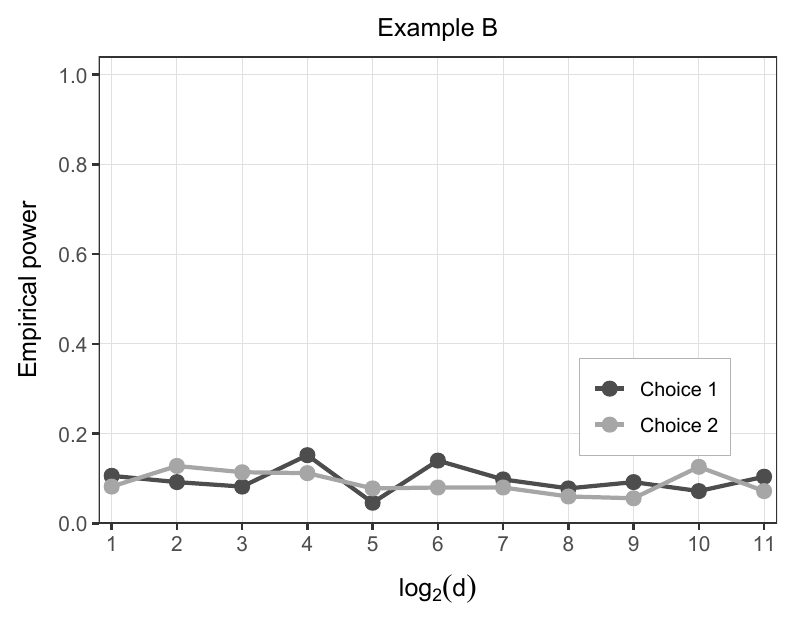}

    \vspace{-0.1in}
    \caption{Powers of the modified rank energy test in Examples A and B for two choices of the reference distribution}
    \label{fig:modified-rank-energy}
\end{figure}

To overcome this problem, one can use a reference distribution with a more heterogeneous structure. For instance, 
a mixture of two normal distributions 
with different locations and scatter matrices can be used. 
Figure 2 shows the performance of the resulting rank energy test for two choices of the reference distribution; Choice-1: an equal mixture of $N_d(5{\bf 1}_d,{\bf I}_d)$ and $N_d(-5{\bf 1}_d,10{\bf I}_d)$ and Choice-2: 
an equal mixture of $N_d(5{\thetavec_d},{\bf I}_d)$ and $N_d(-5{\thetavec}_d,10{\bf I}_d)$, where $\thetavec_d=(\theta_1,\ldots,\theta_d)$ with $\theta_i=(-1)^{i}$ for $i=1,2,\ldots,d$. In Example A, while Choice-1 led to an excellent performance, Choice-2 yielded powers close to the nominal level for all values of $d$. However, none of them could provide any improvement over the rank energy test in Example B.

 The reason can be explained using the distance concentration phenomenon.
For $\Xvec \sim F$ (here $F=N_d(\mathbf{0}_d,{\bf I}_d)$), let the limiting value of
$d^{-1}\|\Xvec-{\bf U}\|^2$ be $\alpha_{1i}$ ($i=1,2$) 
when ${\bf U}\sim R_i$, the
$i$-th component of the reference distribution. For $\Yvec \sim G$ (here $G$ denotes
the location or the scale alternative), let the limiting value of
$d^{-1}\|\Yvec-{\bf U}\|^2$ be $\alpha_{2i}$ for ${\bf U} \sim R_i$ ($i=1,2$). 
If $k$ out of the $m$ observations from $F$ are transported to
$R_1$, the scaled transportation cost converges to
\[
k\alpha_{11}+(m-k)\alpha_{12}+(m_0-k)\alpha_{21}
+(n-m_0+k)\alpha_{22} =m\alpha_{12}+m_0\alpha_{21}+(n-m_0)\alpha_{22}
+k(\alpha_{11}+\alpha_{22}-\alpha_{12}-\alpha_{21}),
\]
where $m_0$ denotes the number of reference observations generated from $R_1$. So, if 
$\alpha_{11}+\alpha_{22}
>
\alpha_{12}+\alpha_{21}$ (respectively, $\alpha_{11}+\alpha_{22}
<
\alpha_{12}+\alpha_{21}$), 
this 
is minimized when $k$ takes its smallest (respectively, largest) possible value. 
So, the optimal transport asymptotically assigns 
almost all observations from $F$ 
to either $R_1$ or $R_2$ and almost all observations from $G$ to the other component, thereby preserving the separation between the two samples. 
However, if the equality $
\alpha_{11}+\alpha_{22}
=
\alpha_{12}+\alpha_{21}
$  holds, all feasible allocations have the same limiting cost. So, one of them is selected randomly. Consequently, the optimal transport may mix the two samples and the rank-energy test loses power. For the location alternative, this equality holds for Choice-2, while for the scale alternative it holds for both choices of the reference distribution. 
This explains the behavior of the rank energy test observed in Figure 2.
         


The performance of the rank energy test in Example B can be improved if instead of $C(\zvec,\wvec)=\|\zvec-\wvec\|^2$, we use $C(\zvec,\wvec)=\|\zvec-\wvec\|$ as the transportation cost. Figure 3 shows the powers of the resulting test in Examples A and B for the two choices of the reference distribution considered before. In Example A, its performance was almost similar to what we observed earlier. But, in Example B, it had much better performance. The power of the test increased with $d$ for both choices of the reference distribution. This can again be explained using the distance concentration phenomenon. 
Proceeding in the same way as before, one can show that the scaled version of the total transportation cost ${\cal C}(\pi)/\sqrt{d}$ is asymptotically minimized at an extreme allocation whenever 
$
\sqrt{\alpha_{11}}
+\sqrt{\alpha_{22}}
\neq
\sqrt{\alpha_{12}}
+\sqrt{\alpha_{21}}.
$ 
Therefore, in such cases, observations from the two high-dimensional distributions are assigned to two different components of the reference distribution. Importantly, \(\alpha_{11}+\alpha_{22}=\alpha_{12}+\alpha_{21}\) does not generally imply the corresponding equality for their square roots, which explains the improved performance in Example B. 


\begin{figure}[t]
    \centering
    \includegraphics[height=1.7in,width=0.4\textwidth]{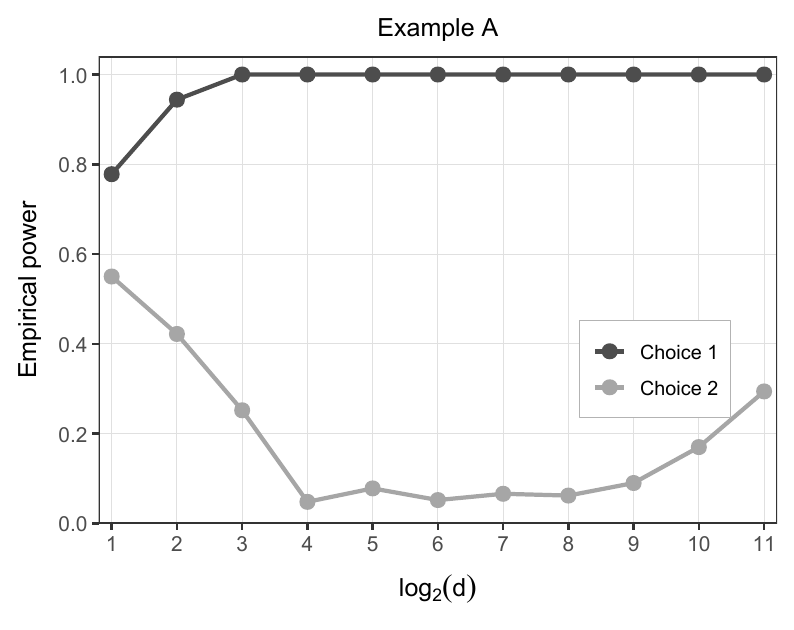}
    \includegraphics[height=1.7in,width=0.4\textwidth]{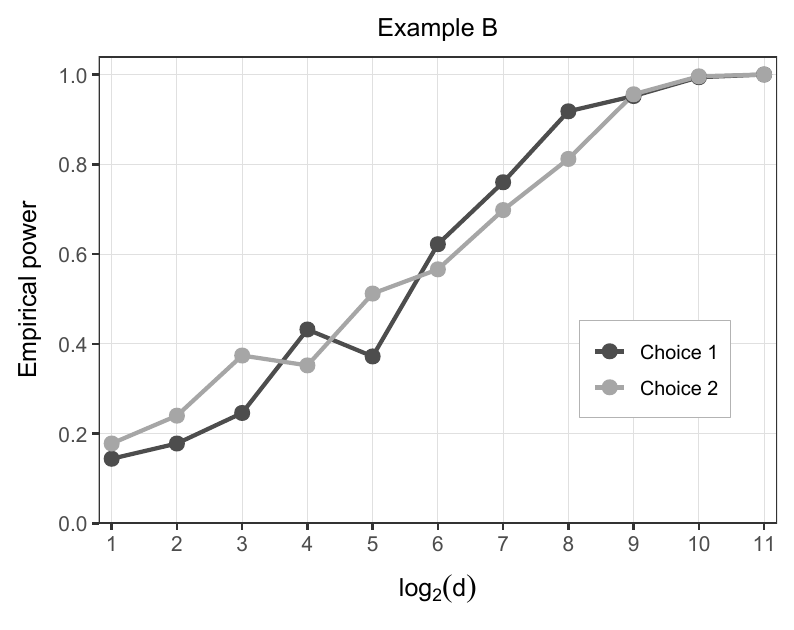}

    \caption{Powers of the modified rank energy test in Examples A and B under the cost function $C(\zvec,\wvec)=\|\zvec-\wvec\|$ for two choices of the reference distributions}\label{fig:modified-rank-energy2}
\end{figure}

\begin{figure}[!b]
    \centering
    \includegraphics[height=1.7in,width=0.4\textwidth]{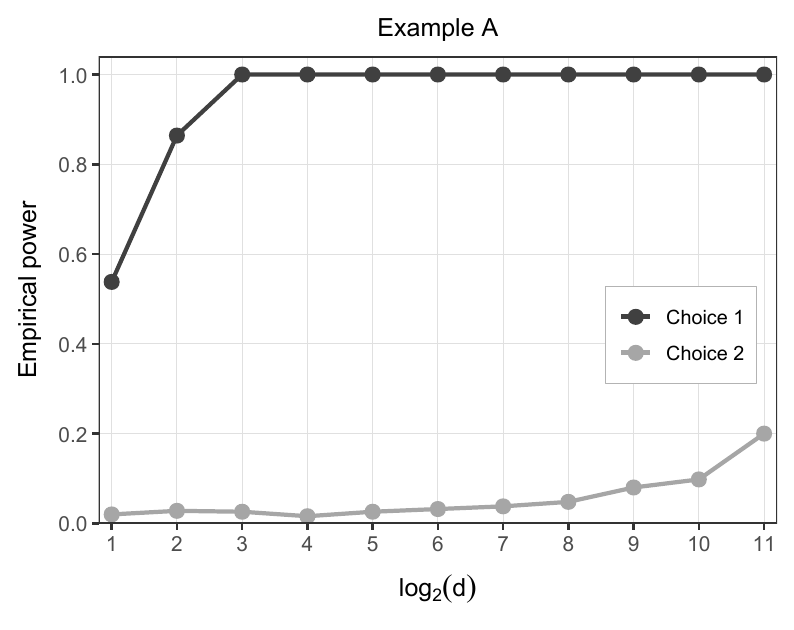}
    \includegraphics[height=1.7in,width=0.4\textwidth]{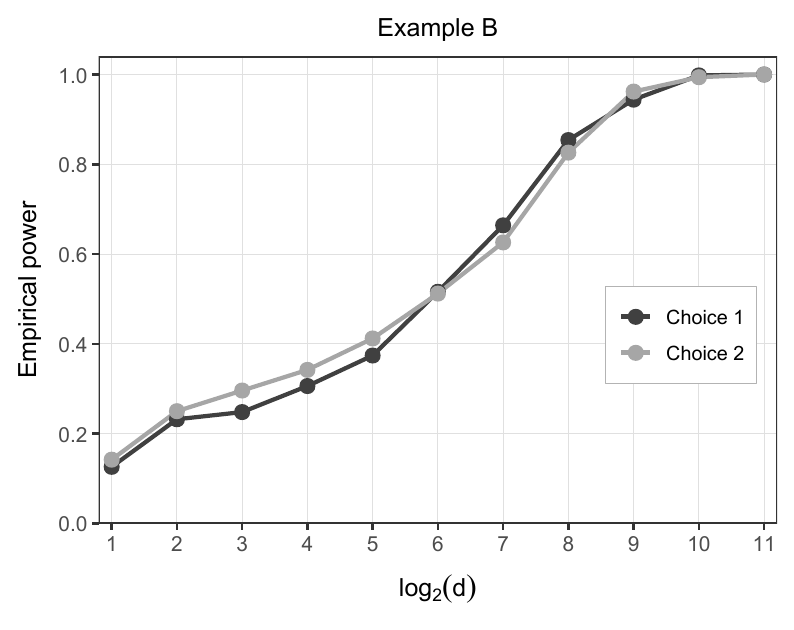}
    
    \caption{Powers of the  test based on TA count in Examples A and B  under the cost function $C(\zvec,\wvec)=\|\zvec-\wvec\|$ for two choices of the reference distributions}
    \label{fig:ta-count}
\end{figure}


The above discussion suggests that the allocation count itself can be used as a test statistic. Let \(k\) denote the number of observations from \(F\) transported to $R_1$. Under the alternative ${\cal H}_1:F\neq G$, \(k\) tends to be unusually large or small, whereas under \({\cal H}_0:F=G\), it has a hypergeometric distribution. We therefore obtain a distribution-free test by rejecting ${\cal H}_0$ for extreme values of \(k\). Figure 4 shows that this test based on transport allocation (TA) count had performance comparable to that of the corresponding rank energy test
in both examples.

The above phenomenon 
can also be observed for distributions with possibly dependent and non-identically distributed coordinate variables if we have at least probability convergence of suitably scaled Euclidean distances. However,
in all such cases, for the success of the rank energy test based on the squared Euclidean distance cost, we need
$\alpha_{11}+\alpha_{22}\neq \alpha_{12}+\alpha_{21}$, whereas for that based on Euclidean distance cost or the corresponding test based on the TA count, we need 
$\sqrt{\alpha_{11}}+\sqrt{\alpha_{22}}\neq \sqrt{\alpha_{12}}+\sqrt{\alpha_{21}}$. Thus, the high-dimensional performance of these tests depends critically on the choice of both the reference distribution and the transportation cost. Since these quantities depend on the unknown underlying distributions, choosing a suitable reference distribution is generally difficult. This motivates our reference-free approach based on optimal derangements in the next section.

\section{Tests based on optimal derangement}


For the moment, suppose that additional observations from $F$ and $G$ are available, and we can use $R_1=F$ and $R_2=G$ as the two components of the reference distribution. Let $\muvec_1$ and $\sigmat_1$ be the mean vector and the dispersion matrix of $F$, while $\muvec_2$ and $\sigmat_2$ are those of $G$. Also assume that $\nu^2$, $\sigma_1^2$ and $\sigma_2^2$ are the limiting (as $d$ tends to infinity) values of $d^{-1}\|\muvec_1-\muvec_2\|^2$, $d^{-1}\text{trace}(\sigmat_1)$ and $d^{-1}\text{trace}(\sigmat_2)$, respectively. Then, under the assumption of convergence of scaled Euclidean distances, one gets $\alpha_{11}=2\sigma_1^2$, $\alpha_{22}=2\sigma_2^2$ and $\alpha_{12}=\alpha_{21}=\sigma_1^2+\sigma_2^2+\nu^2$. Hence, we have $$\sqrt{\alpha_{11}}+\sqrt{\alpha_{22}} =(\sigma_1+\sigma_2)\sqrt{2} \mbox{     and      } \sqrt{\alpha_{12}}+\sqrt{\alpha_{21}}=2\sqrt{\nu^2+\sigma_1^2+\sigma_2^2}.$$ 
Note that the identity  $4({\nu^2+\sigma_1^2+\sigma_2^2})=4\nu^2 + 2(\sigma_1+\sigma_2)^2+ 2(\sigma_1-\sigma_2)^2$
implies $\sqrt{\alpha_{12}} + \sqrt{\alpha_{21}} \ge \sqrt{\alpha_{11}} + \sqrt{\alpha_{22}}$, where the equality holds if and only if $\nu^2=0$ and $\sigma_1^2=\sigma_2^2$. 
So, if the average of the coordinate-wise location difference or that of the scale difference between $F$ and $G$ is asymptotically non-negligible, one can expect the TA count for the optimal transport to be large, indicating a deviation from ${\cal H}_0$. 

Since \(F\) and \(G\) are unknown, however, they cannot be used to construct the reference distribution. An alternative is to randomly split the pooled sample into two halves and use one half as the reference sample. However, this introduces randomness into the procedure and often leads to a substantial loss of power. We therefore avoid sample splitting and instead transport the pooled sample $\{\zvec_1,\zvec_2,\ldots,\zvec_{m+n}\}$ to itself. But, in this case, the identity mapping would trivially minimize the transportation cost. So we restrict attention to permutations \(\pi\) satisfying \(\pi(i)\ne i\) for all \(i=1,2,\ldots,m+n\) and choose the one 
that minimizes the total transportation cost $\sum_{i=1}^{m+n}\|\zvec_i-\zvec_{\pi(i)}\|$. Such a permutation or arrangement is referred to as an optimal derangement. 
  The optimization is a modified assignment problem, with infinite cost assigned to self-matches, and can be solved in polynomial time using the Hungarian algorithm \citep{kuhn1956variants}.

Let \(k\) denote the number of observations from \(F\) that are transported to another observation from \(F\). Then the number of cross-matches is 
$D_{m,n}=2(m-k)$. Under \({\cal H}_1:F\neq G\), we expect \(k\) to be large, or equivalently \(D_{m,n}\) to be small, reflecting the tendency of the optimal derangement to connect observations from the same distribution. We therefore reject \({\cal H}_0:F=G\) for large \(k\) (equivalently, small \(D_{m,n}\)). Unlike the rank energy test, this procedure does not require an external reference distribution. We call it the Optimal DERangement (ODER) test.

Note that an optimal derangement may consist of cycles of different sizes. For example, with six observations, possible cycle structures include three cycles of size 2, two cycles of size 3, two cycles of sizes 2 and 4, and a single cycle of size 6; these are illustrated in Figure 5.
If we restrict to derangements with cycles of size $2$ only, finding optimal derangements is equivalent to finding the pairs of observations with minimum total distance. In that case, $D_{m,n}$ turns 
out to be twice the number of cross-match pairs, and the test coincides with the distribution-free test of \cite{MinDist}. On the other hand, if we restrict derangements to a single cycle of size $m+n$, the optimization problem reduces to finding a shortest Hamiltonian cycle (SHC), and \(D_{m,n}\) gives the number of runs in that cycle. This construction is closely related to the distribution-free runs test of Biswas et al. (2014), although the two tests are not identical.

\begin{figure}[h]
    \centering
    \includegraphics[height=1.5in,width=0.9\textwidth]{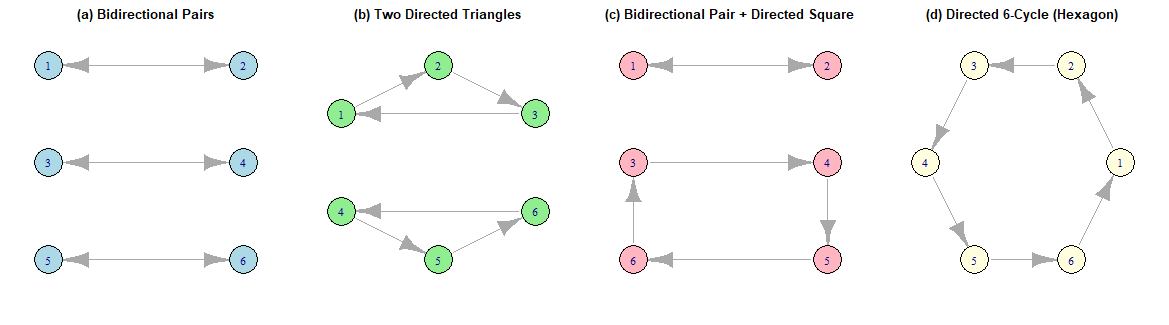}
    \vspace{-0.25in}
    \caption{Different types of derangements on a set of six observations}
    \vspace{-0.1in}
    \label{fig:derangement-illustration}
\end{figure}

If we fix the size of the derangement cycles or fix the cycle structure of the derangements, exchangeability of $\zvec_1,\zvec_2,\ldots,\zvec_{m+n}$ under ${\cal H}_0$ implies that all derangements within the specified class are equally likely, and hence the resulting test is distribution-free. 
For instance, when $N=m+n$ is even, restriction to cycles of size 2 gives the distribution-free MDP test of Rosenbaum (2005), which can be computed in polynomial time \citep{lu2011optimal}. But for cycles of size $\ell>2$ or any other fixed cycle structures, the corresponding optimization is generally NP-hard. However, for derangements consisting of a single cycle of size \(m+n\), the optimization problem reduces to finding an SHC. We can approximate this using the cheapest-link algorithm 
\citep{tannenbaum2010excursions} or the Christofides algorithm \citep{christofides2022}. Henceforth, we shall refer to this test as the ODER-SHC test, and the corresponding test statistic will be denoted by $D_{m,n}^{\circ}$. Since
two possible orientations of a Hamiltonian cycle have the same cost and yield the same test statistic, in practice, one of them can be chosen randomly.

However, we do not face such computational issues for the ODER test, where no restriction is imposed on the cycle structure. As we have mentioned before, in such cases, the optimal derangement can be found in polynomial time using the Hungarian algorithm   \cite{kuhn1956variants}. The optimal derangement may contain one or more cycles of size bigger than 2. Since both orientations of such a cycle have the same cost and yield the same number of cross-matches, we can work with either of them. However, unlike the tests based on a fixed cycle structure, the ODER test is not distribution-free. Although derangements within each cycle type are equally likely under \({\cal H}_0\), the probability of each cycle type can depend on the common underlying distribution. We therefore use a permutation method to determine the critical value. Importantly, for a fixed pooled sample, the optimal derangement does not change under permutations of the class labels, so the optimization needs to be performed only once. This leads to a substantial reduction in computational cost.


Figure 6 shows that both ODER and ODER-SHC tests had excellent performance in Examples A and B, with ODER-SHC exhibiting higher power than ODER, particularly in Example B. 
\begin{figure}[htbp]
    \centering

    \includegraphics[height=1.7in,width=0.4\textwidth]{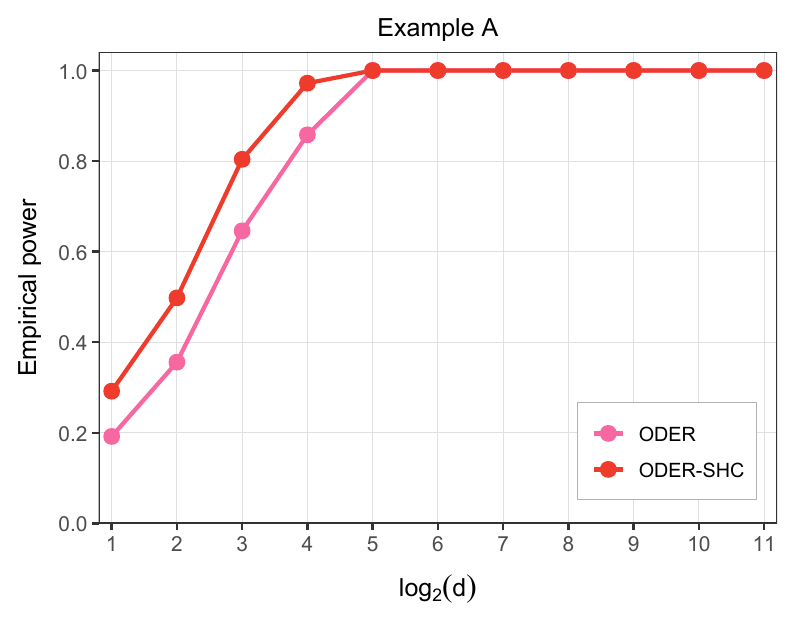}
    \includegraphics[height=1.7in,width=0.4\textwidth]{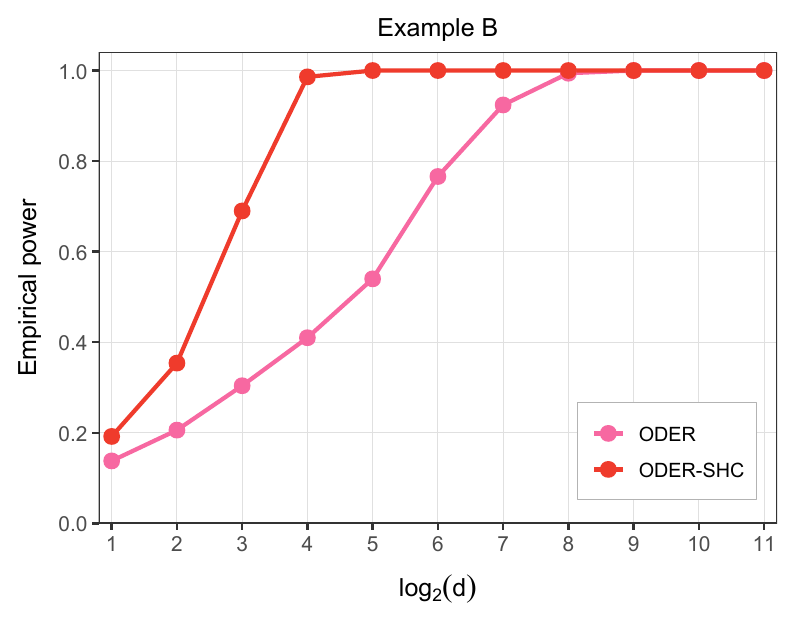}

    \vspace{-0.05in}

    \caption{Powers of ODER and ODER-SHC tests in Examples A and B}
    \vspace{-0.1in}
    \label{fig:oder-ab}
\end{figure}



\subsection {Behavior of the proposed tests for high dimension, low sample size data}
\vspace{0.1in}

Let $\Xvec_1,\Xvec_2,\ldots,\Xvec_m$ be independent copies of $\Xvec=(X^{(1)},\ldots, X^{(d)})\sim F$ and $\Yvec_1,\Yvec_2,\ldots,\Yvec_n$ be independent copies of $\Yvec=(Y^{(1)},\ldots, Y^{(d)})\sim G$. To investigate the high-dimensional behavior of the derangement-based tests, we consider the following assumptions. 
\begin{enumerate}
    \item[(A1)] Fourth moments of the coordinate variables under $F$ and $G$ are uniformly bounded.
    \item[(A2)] For ${\bf W}=\Xvec_1-\Xvec_2$, $\Yvec_1-\Yvec_2$ and $\Xvec_1-\Yvec_1$,
    $\sum_{q}\sum_{r} \Big|corr\big((W^{(q)})^2,(W^{(r)})^2\big)\Big|$ is of the order $o(d^2)$.
    \item[(A3)] There exist $\nu^2$, $\sigma_1^2$ and $\sigma^2_{2}$ such that $~d^{-1}\|\muvec_1-\muvec_2\|^2 \rightarrow \nu^2$, $d^{-1}{\text tr}({\sigmat_1}) \rightarrow \sigma_1^2$ and $d^{-1}{\text tr}({\sigmat_2}) \rightarrow \sigma_2^2$ as $d \rightarrow \infty$. 
\end{enumerate}
Under Assumptions (A1) and (A2), the weak law of large numbers (WLLN) holds for the sequence of
random variables $\{(W^{(q)})^2:q\ge 1\}$ (i.e., $\Big|{d}^{-1}\|{\bf W}\|^2-E\big({d}^{-1}\|{\bf W}\|^2\big)\Big| \stackrel{P}{\rightarrow}0$ as $d \rightarrow \infty$), and hence we have convergence of suitably scaled pairwise distances. 
These assumptions are quite common in the HDLSS literature. \cite{Geometry} considered the $d$-dimensional observations as time series truncated at time $d$. In addition to (A1), they assumed the $\rho$-mixing property of the time series to study the high-dimensional behavior of some nonparametric classifier. The weak dependence assumption (A2) holds under that $\rho$-mixing condition. Other similar conditions for weak dependence and convergence of pairwise distances can be found in \cite{ahn2007high,jung2009pca,aoshima2018two} and \cite{banerjee2025high}. 
However, in the case of independent and identically distributed (i.i.d.) coordinate variables (as in Examples A and B), WLLN holds if $F$ and $G$ have finite second moments; Assumptions (A1) and (A2) are not needed there. Assumption (A3) gives the limiting value of $\E\big(d^{-1}\|{\bf W}\|^2\big)$ and hence that of 
$d^{-1}\|{\bf W}\|^2$.  For ${\bf W}=\Xvec_1-\Xvec_2$, $\Yvec_1-\Yvec_2$ and $\Xvec_1-\Yvec_1$, these limiting values are $2\sigma_1^2$, $2\sigma_2^2$ and $\nu^2+\sigma_1^2+\sigma_2^2$, respectively.
Note that if $\sigma_1^2 \neq \sigma_2^2$ or $\nu^2>0$, we have $\Delta=2\sqrt{\nu^2+\sigma_1^2+\sigma_2^2}-\sigma_1\sqrt{2}-\sigma_2\sqrt{2}>0$. So, under (A1)–(A3), each additional pair of cross-sample edges increases the cost of a derangement relative to one with fewer cross-sample edges. Consequently, 
$D_{m,n}$ and $D_{m,n}^{\circ}$ attain their respective minimum values $0$ and $2$, respectively,  with probability tending to $1$ as $d$ grows to infinity (see the proof of Theorem 3.1). So, if the probabilities of these extreme events
under ${\cal H}_0$ are smaller than the level of the tests $\alpha$, the powers of these tests will converge to unity. Since the ODER-SHC test is distribution-free, it is easy to find this probability. But for the ODER test based on permutation distribution, this probability depends on the cycle structure under optimal derangement. However, it is possible to find a uniform upper bound for this probability, which is given by $B_{m,n}=\big.{
\sum_{j=1}^{\lfloor r/2\rfloor}
\binom{\lfloor N/2\rfloor}{j}
}\big/
\binom{N}{r}$, where $r=\min\{m,n\}$.
Consequently, we have the following theorem asserting 
the high-dimensional consistency of the derangement-based tests. 
Since finding the SHC is NP-hard, we approximate the SHC using the cheapest-link algorithm in all numerical work. All theoretical results involving the ODER-SHC test are stated initially for the exact SHC. However, as shown in Remarks 1 and 2 in the Appendix, the corresponding consistency results also hold for the implemented version of the test based on the cheapest-link approximation.


\begin{thm}
\label{thm:HDLSS-derangement}
Suppose $\Xvec_1,\Xvec_2,\ldots,\Xvec_m$ ($m\ge 2$) are independent copies of $\Xvec \sim F$, and $\Yvec_1,\Yvec_2,\ldots,\Yvec_n$ ($n \ge 2$) are independent copies of $\Yvec \sim G$, where $F$ and $G$ satisfy (A1)-(A3) with $\nu^2>0$ or $\sigma_1^2\neq \sigma_2^2$.
\begin{enumerate}
\item[(a)] 
If $B_{m,n}<\alpha$, the power of the ODER test at level $\alpha$ converges to $1$ as $d$ tends to infinity.
\item[(b)] If $m!n!/(m+n-1)!<\alpha$, the power of the ODER-SHC test at level $\alpha$ converges to 1 as $d$ tends to infinity. 
\end{enumerate}
\end{thm}

One can show that $B_{m,n}<0.05$ for all $m=n\ge 6$. So, for the high-dimensional consistency of the ODER test at $5\%$ level, it is enough to have $6$ observations from each distribution. Since $m!n!/(m+n-1)!<0.05$ for all $m=n\ge 5$, for the ODER-SHC test, this consistency holds for all $m=n\ge 5$. As we have mentioned before, the consistency result stated in Theorem 3.1(b) continues to hold for the implemented version of the test based on the approximate SHC obtained using the cheapest-link algorithm (see Remark 1 in the Appendix after the proof of Theorem 3.1). 



Assuming the independence and normality of the coordinate variables, \cite{BorysovHannigMarron2014} established distance convergence results allowing the pooled sample size $N=m+n$ to increase with the dimension 
$d$ at a sub-exponential rate. Under sub-Gaussianity and $m$-dependence of the coordinate variables (under some permutation), we can prove similar distance convergence results and hence the high-dimensional consistency of ODER and ODER-SHC tests as in Theorem 3.1 for $\log N=o(d)$ (see Lemma A.1 in the Appendix).

The condition $\nu^2>0$ or $|\sigma_1^2-\sigma_2^2|>0$ in Theorem 3.1 requires the location difference $\|\muvec_{1}-\muvec_{2}\|^2$ or the scale difference $\left|{\text tr}(\sigmat_{1})-{\text tr}(\sigmat_{2})\right|$ to grow linearly with $d$. This may not hold if $F$ and $G$ differ only in $d^{\beta}$ coordinates with $\beta<1$. However, our proposed tests may have high-dimensional consistency
even for such weaker or potentially sparse alternatives. For four independent random vectors $\Xvec_1,\Xvec_2{\sim}F$ and $\Yvec_1,\Yvec_2 {\sim} G$,
define $S_{FF,d}^2$, $S_{GG,d}^2$ and $S_{FG,d}^2$ as the expectations of $\|\Xvec_1-\Xvec_2\|^2/d$, $\|\Yvec_1-\Yvec_2\|^2/d$ and $\|\Xvec_1-\Yvec_1\|^2/d$, respectively. Let $V^2_{FF,d}$, $V^2_{GG,d}$ and $V^2_{FG,d}$ be their respective variances. While $S_d=2S_{FG,d}-S_{FF,d}-S_{GG,d}$ gives a signal against ${\cal H}_0$,  $V_d=\max\{V_{FF,d},V_{GG,d},V_{FG,d}\}$ can be viewed as the associated stochastic variation or noise.
For the high-dimensional consistency of the proposed tests, we only need $S_d$ to be of higher asymptotic order than $V_d$ (i.e., $S_d \gg V_d$). This result is stated below.

\begin{thm}
\label{thm:joint-consistency}
    Suppose $\Xvec_1,\Xvec_2,\ldots,\Xvec_m(m\geq 2)$ are independent copies of $\Xvec \sim F$, and $\Yvec_1,\Yvec_2,\ldots,\Yvec_n(n\geq 2)$ are independent copies of $\Yvec \sim G$, where $F$ and $G$ satisfy (A1)-(A3).   
    If $\lim_{d \rightarrow \infty} d^{-1}{\text tr}(\sigmat_i)>0$ for $i=1,2$, and $S_d/V_d \rightarrow \infty$ as $d \rightarrow \infty$, then we have the following results.
\begin{enumerate}
\item[(a)]  
If $B_{m,n}<\alpha$, the power of the ODER test of level $\alpha$ converges to $1$ as $d$ tends to infinity.
\item[(b)] If $m!n!/(m+n-1)!<\alpha$, the power of the ODER-SHC test of level $\alpha$ converges to 1 as $d$ tends to infinity.  
\end{enumerate}
\end{thm}

Theorem 3.2(b) holds even when the SHC is approximated using the cheapest-link algorithm (see Remark 2 in the Appendix after the proof of Theorem 3.2). Note that if the coordinate variables in \(F\) and \(G\) are independent or \(m\)-dependent under some permutation, $V_d$ is of the order $d^{-1/2}$ (i.e., $V_d \asymp d^{-1/2}$). In such cases, instead of requiring \(S_d\) to converge to a positive constant (as needed in Theorem 3.1), it is enough to have the divergence of 
\(d^{1/2}S_d\). Thus, the result holds even when
 $\|\muvec_{1}-\muvec_{2}\|^2 \gg d^{1/2}$ or $|\text{tr}(\sigmat_{1})-\text{tr}(\sigmat_{2})|\gg d^{3/4}$ (see Lemma A.2 in the Appendix), and we do not need them to grow linearly with $d$.

Now, one may ask what happens to these tests when the sufficient conditions for consistency in Theorems 3.1 and 3.2 are not satisfied. To investigate this, we consider the following two examples, both with $\nu^2=0$ and $\sigma_1^2 = \sigma_2^2$. 

{\bf Example C:} Each of the two distributions has i.i.d. coordinate variables. For the first distribution, they are i.i.d. \(N(0,2)\), whereas for the second, they are i.i.d. \(t_4\), the standard \(t\)-distribution with 4 degrees of freedom.

{\bf Example D:} As in Example B, consider two $d$-variate normal distributions with the same mean vector ${\bf 0}_d$ but different dispersion matrices ${\sigmat}_{1d}$ and ${\sigmat}_{2d}$.
The matrix \({\sigmat}_{1d}\) is diagonal, with its first \(d/2\) diagonal elements equal to \(1\) and its remaining diagonal elements equal to \(2\), whereas \({\sigmat}_{2d}\) is diagonal, with its first \(d/2\) diagonal elements equal to \(2\) and its remaining diagonal elements equal to \(1\).

\begin{figure}[h]
    \centering
    \includegraphics[height=1.7in,width=0.4\textwidth]{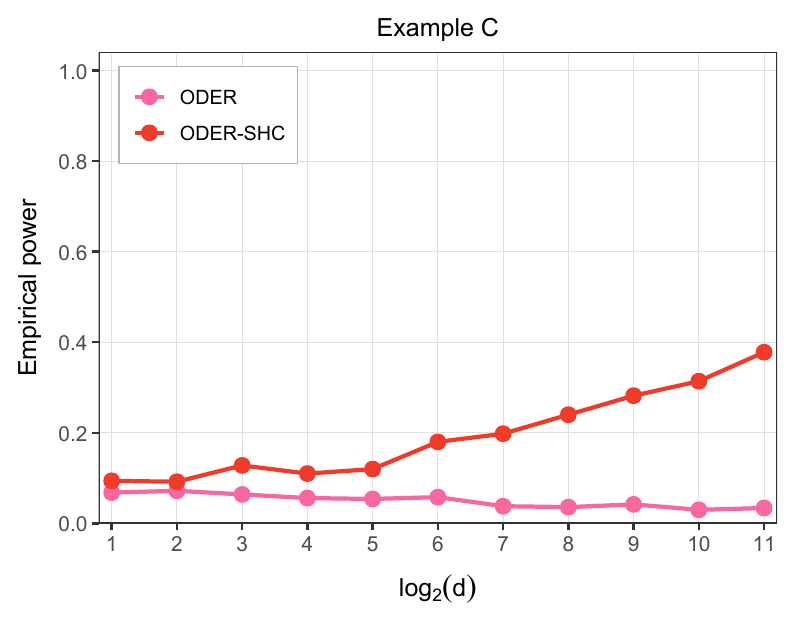}
    \includegraphics[height=1.7in,width=0.4\textwidth]{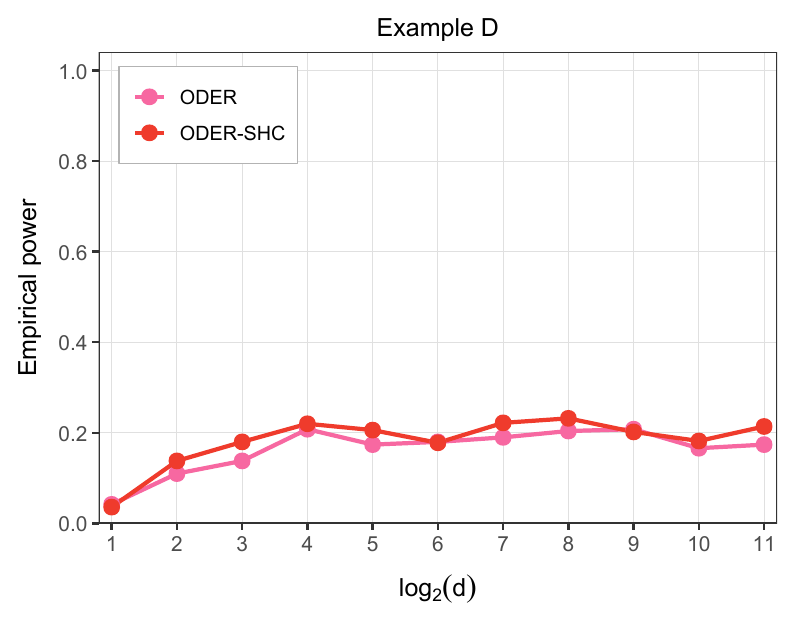}

    \vspace{-0.05in}
    \caption{Powers of ODER and ODER-SHC tests in Examples C and D}
    \label{fig:oder-cd}
\end{figure}

In these two examples, since the marginal distribution of each coordinate variable differs under the two competing distributions, one may expect the power of a test to increase with the dimension. But, in Example D, the powers of ODER and ODER-SHC tests remained unsatisfactory as $d$ increased (see Figure 7). In Example C, the power of the ODER-SHC test showed a slowly increasing trend with \(d\), but the power of the ODER test remained close to the nominal level. These observations suggest that the choice of the distance function may play an important role when the conditions of Theorems 3.1 and 3.2 are not satisfied. This motivates us to investigate the \(\ell_1\)-distance as an alternative to the Euclidean distance in the next subsection.

\subsection{ODER and ODER-SHC tests based on $\ell_1$-distance}
\vspace{0.1in}

When we used the $\ell_1$-distance as the transportation cost (i.e., $C(\zvec,\wvec)=\|\zvec-\wvec\|_1$), the corresponding ODER and ODER-SHC tests (henceforth referred to as ODER-L1 and ODER-SHC-L1 tests, respectively), both exhibited excellent performance in Examples C and D. Figure 8 shows the powers of these tests based on $\ell_1$-distance in Examples A-D along with those of the corresponding tests based on Euclidean distance reported earlier. This figure clearly shows that the $\ell_1$-distance based tests also performed well in Examples A and B. 

\begin{figure}[htbp]
    \centering

    \begin{tabular}{cc}
        \includegraphics[height=1.6in, width=0.4\textwidth]{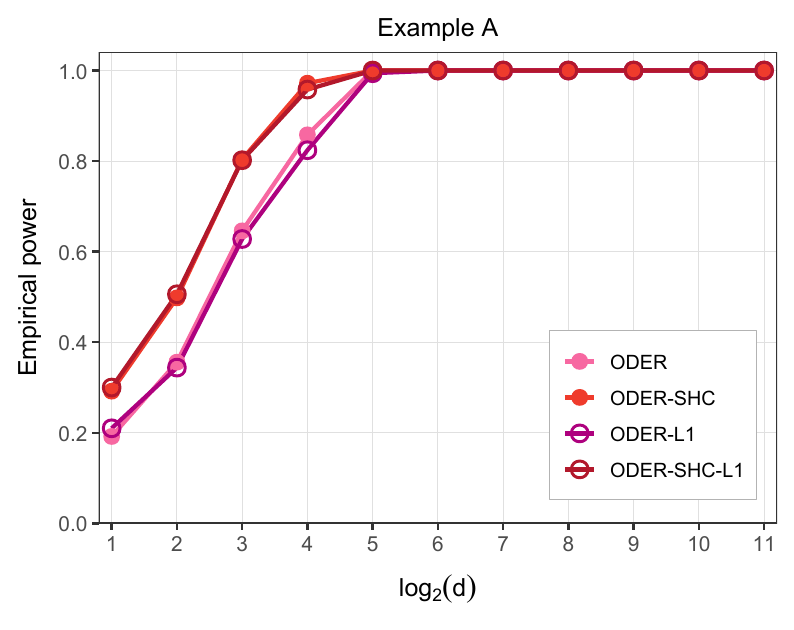} &
        \includegraphics[height=1.6in, width=0.4\textwidth]{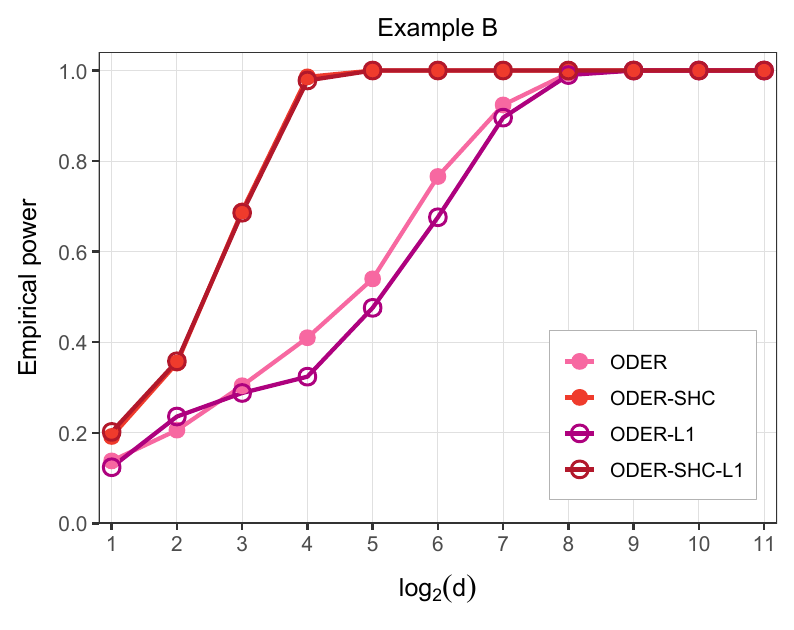} \\
        
    \includegraphics[height=1.6in, width=0.4\textwidth]{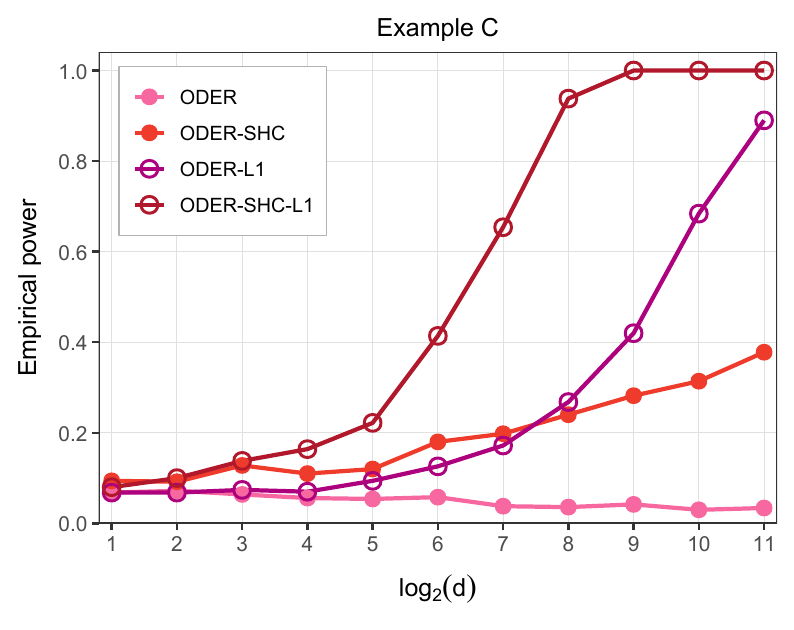} &
        \includegraphics[height=1.6in,width=0.4\textwidth]{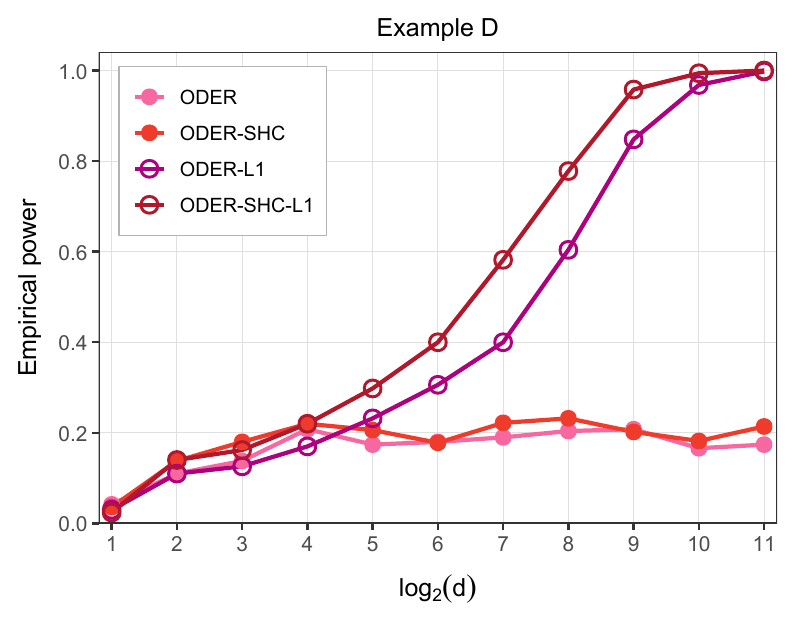}
\end{tabular}



    \vspace{-0.1in}

    \caption{
    Empirical powers of derangement-based tests 
    in Examples A, B, C and D.
    }
    \vspace{-0.2in}
    \label{fig:examples-ABCD-all4}
\end{figure}

This motivates us to investigate the high-dimensional behavior of ODER-L1 and ODER-SHC-L1 tests. For this purpose, we consider the following assumptions analogous to (A1)-(A3).

\begin{enumerate}
    \item[(A1$^{\circ}$)] Second moments of the coordinate variables under $F$ and $G$ are uniformly bounded.
    
    \item[(A2$^{\circ}$)] For ${\bf W}=\Xvec_1-\Xvec_2$, $\Yvec_1-\Yvec_2$ and $\Xvec_1-\Yvec_1$,
    $\sum_{q}\sum_{r} \Big|corr\big(|W^{(q)}|,|W^{(r)}|\big)\Big|$ is of the order $o(d^2)$.

    \item[(A3$^{\circ}$)] $d^{-1}E\|{\bf W}\|_1$ has finite limits for ${\bf W}=\Xvec_1-\Xvec_2$, $\Yvec_1-\Yvec_2$ and $\Xvec_1-\Yvec_1$, 
\vspace{-0.1in}
\end{enumerate}

These assumptions ensure WLLN for the sequence $\{|W^{(q)}|:~q\ge 1\}$, and hence
${d}^{-1}\|{\bf W}\|_1-{d}^{-1}E\|{\bf W}\|_1 \stackrel{P}{\rightarrow}0$ as $d\rightarrow \infty$. Now, for each $q=1,2,\ldots,d$, we have
${\cal E}_{F,G}^{(q)}=2E|X_1^{(q)}-Y_1^{(q)}|-E|X_1^{(q)}-X_2^{(q)}|-E|Y_1^{(q)}-Y_2^{(q)}|\ge 0$, where the equality holds if and only if the $q$-th coordinate variables in $F$ and $G$ have the same distribution \citep[follows from][]{BARINGHAUS2004}. So, for any finite $d$, ${d}^{-1}\sum_{q=1}^{d}{\cal E}_{F,G}^{(q)}>0$ unless $F$ and $G$ have the same set of coordinate-wise marginals. Therefore, under ${\cal H}_1$, it is somewhat reasonable to assume that $\lim_{d \rightarrow \infty}d^{-1}\sum_{q=1}^{d}{\cal E}_{F,G}^{(q)}>0$. Here ${\cal E}_{F,G}^{(q)}$ is the energy distance \citep[see, e.g.,][]{szekely2013energy} between the $q$-th marginals of $F$ and $G$. So, it assumes that the average of the coordinate-wise energy distances between $F$ and $G$ is asymptotically non-negligible. Under this assumption, we have the high-dimensional consistency of ODER-L1 and ODER-SHC-L1 tests.

\begin{thm}
\label{thm:HDLSS-L1}
Suppose $\Xvec_1,\Xvec_2,\ldots,\Xvec_m(m\geq 2)$ are independent copies
of $\Xvec\sim F$, and
$\Yvec_1,\Yvec_2,\ldots,\Yvec_n(n\geq 2)$ are independent copies of
$\Yvec\sim G$, where $F$ and $G$ satisfy Assumptions (A1$^{\circ}$)- (A3$^{\circ}$) with $\lim_{d \rightarrow \infty}\frac{1}{d}\sum_{q=1}^{d}{\cal E}_{F,G}^{(q)}>0$. Then, we have the following results for ODER-L1 and ODER-SHC-L1 tests at level $\alpha$ $(0<\alpha<1)$.
\begin{enumerate}
\item[(a)]  
If $B_{m,n}<\alpha$, the power of the ODER-L1 test converges to $1$ as $d$ tends to infinity.
\item[(b)] If $m!n!/(m+n-1)!<\alpha$, the power of the ODER-SHC-L1 test converges to 1 as $d$ tends to infinity. 
\end{enumerate}
\end{thm}

Note that Examples C and D satisfy the condition of Theorem 3.3, which explains the superior performance of ODER-L1 and ODER-SHC-L1 tests. As in the case of Theorem 3.1, one can show that if the coordinate variables in $F$ and $G$ are sub-Gaussian and $m$-dependent under some permutation, Theorem 3.3 holds even when $N$ increases with $d$ at a sub-exponential rate (i.e., $\log N=o(d)$). Since the proof is essentially the same as that of Lemma A.1, we omit it.

Theorem 3.3 assumes that $\sum_{q=1}^{d} {\cal E}_{F,G}^{(q)}$, the total coordinate-wise energy distance between $F$ and $G$, increases at least linearly with $d$. However, this sufficient condition can be relaxed
in some situations. Let $V^{\circ2}_{FF,d}$, $V^{\circ2}_{GG,d}$ and $V^{\circ2}_{FG,d}$ denote the variances of $\|\Xvec_1-\Xvec_2\|_1/d$, $\|\Yvec_1-\Yvec_2\|_1/d$ and $\|\Xvec_1-\Yvec_1\|_1/d$, respectively, and let $V_d^{\circ}=\max\{V_{FF,d}^{\circ},V_{GG,d}^{\circ},V_{FG,d}^{\circ}\}$.  Then, for the high-dimensional consistency of ODER-L1 and ODER-SHC-L1 tests, it is sufficient that 
 $S_d^{\circ}=\frac{1}{d}\sum_{q=1}^{d} {\cal E}_{F,G}^{(q)}\gg V_d^{\circ}$  (see Theorem 3.4 below). So, if the coordinate variables in $F$ and $G$ are independent or $m$-dependent (under some permutation), in which case $V_d^{\circ}$ is of the order $O(d^{-1/2})$, unlike Theorem 3.3, $S_d^{\circ}$ need not converge to a positive constant; it is sufficient that $d^{1/2}S_d^{\circ} \rightarrow \infty$ as $d \rightarrow \infty$. 

\begin{thm}
\label{thm:HDLSS-L1-sparse}
    Suppose $\Xvec_1,\Xvec_2,\ldots,\Xvec_m(m\geq 2)$ are independent copies of $\Xvec \sim F$, and $\Yvec_1,\Yvec_2,\ldots,\Yvec_n(n\geq 2)$ are independent copies of $\Yvec \sim G$, where $F$ and $G$ satisfy Assumptions (A1$^{\circ}$)-($A3^{\circ}$). If $S_d^{\circ}/V_d^{\circ} \rightarrow \infty$ as $d \rightarrow \infty$, then we have the following results for ODER-L1 and ODER-SHC-L1 tests at level $\alpha$ $(0<\alpha<1)$.
\begin{enumerate}
\item[(a)]  
If $B_{m,n}<\alpha$, the power of the ODER-L1 test converges to $1$ as $d$ tends to infinity.
\item[(b)] If $m!n!/(m+n-1)!<\alpha$, the power of the ODER-SHC-L1 test converges to 1 as $d$ tends to infinity. 
\end{enumerate}
\end{thm}

As before, Theorem 3.3(b) and Theorem 3.4(b) are stated for the test based on the true SHC. However, the HDLSS consistency of the ODER-SHC-L1 test continues to hold for the implemented version of the test based on the cheapest-link algorithm. This is discussed in Remarks 3 and 4 in the Appendix following the proofs of these theorems.

\section{Results from the analysis of benchmark datasets}

We analyze some benchmark datasets to compare the empirical performance of the derangement-based tests with some existing two-sample tests. In particular, we consider the test based on the Energy statistic \citep{BARINGHAUS2004}, the kernel two-sample test based on Maximum Mean Discrepancy \citep{Kernel}, the multivariate run test based on the Minimum Spanning Tree \citep{MST}, the test based on Nearest Neighbors \citep{Neighbors1986,Neighbors1988}, and the test based on Minimum Distance Pairing \citep{MinDist} for comparison. 
Henceforth, we refer to them as the Energy, MMD, MST-run, NN and MDP tests, respectively. The MDP test has the distribution-free property. In all other cases, we use the conditional test based on the permutation principle, with critical values computed from $500$ random permutations. The codes for these competing tests are available in different {\tt R} packages. For the derangement-based tests, we use our own {\tt R} codes, which are available from the first author upon request.
For the MMD test, we use the radial basis function kernel, and the bandwidth is selected using the median heuristic \citep[see][]{Kernel}. For the nearest neighbor method, following \cite{Geometry}, we use $k=1$ for these high-dimensional data. Throughout this article, for all numerical work, all tests are considered at 5\% nominal level. 

In this section, we analyze six datasets. Among them, Trace data and Lightning-2 data are available at the UCR Time Series Classification Archive. 
Colon data can be obtained in various {\tt R} packages, including {\tt plsgenomics}. The other datasets, Ionosphere data, Sonar data and Arcene data, are taken from the UCI Machine Learning repository. Brief descriptions of these datasets are given below.

{\bf Ionosphere data:} This radar dataset, collected in Goose Bay, Labrador, contains 351 observations, each with 34 attributes derived from radar measurements. The observations comprise 225 `good' radar returns, indicating signals showing evidence of structure in the ionosphere, and 126 `bad' radar returns, indicating 
signals without such structure.

{\bf Sonar data:} This dataset consists of 208 sonar observations, with 111 signals reflected from metal cylinders (mines) and 97 from rocks. Each observation is represented by 60 values between 0 and 1 corresponding to energy measurements in different frequency bands.

{\bf Trace data:} This synthetically generated dataset simulates transient instrumentation failures in a nuclear power plant. It contains four classes, each with 50 time-series observations of dimension 275. For our analysis, we consider observations from the first two classes.

{\bf Lightning-2 data:} This dataset contains observations from two classes of lightning events: cloud-to-ground and intra-cloud lightning. Each observation is represented by 637 features obtained from the recorded electromagnetic signal. The dataset contains 48 observations from one class and 73 from the other.

{\bf Colon data:} Originally collected by \cite{alon1999broad}, this dataset contains gene-expression measurements for 62 colon tissue samples, with 2,000 gene-expression variables. Of the 62 samples, 40 are from tumor tissues and 22 from normal tissues. We treat these two groups as the two samples.

{\bf Arcene data:} This dataset was developed for the NIPS 2003 feature selection challenge by combining three mass-spectrometry datasets. Each observation has 10,000 numerical features, including 7,000 mass-spectrometry features and 3,000 synthetic distractors. The dataset contains 900 observations, with 398 and 502 observations in the two classes, respectively.

\begin{figure}[ht]
    \centering


    \begin{minipage}[t]{0.31\textwidth}
        \centering
        \includegraphics[
            width=\linewidth,
            height=1.7in,
            keepaspectratio
        ]{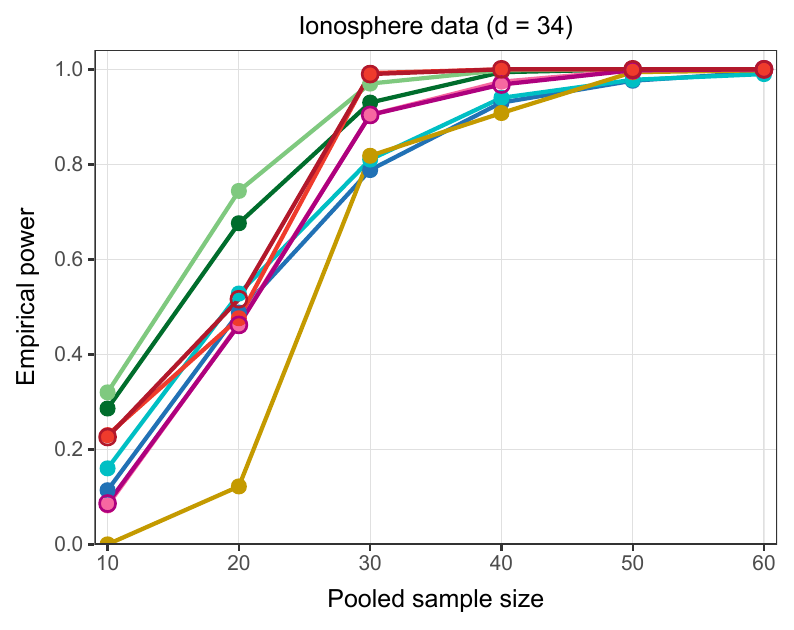}
    \end{minipage}
    \hfill
    \begin{minipage}[t]{0.31\textwidth}
        \centering
        \includegraphics[
            width=\linewidth,
            height=1.7in,
            keepaspectratio
        ]{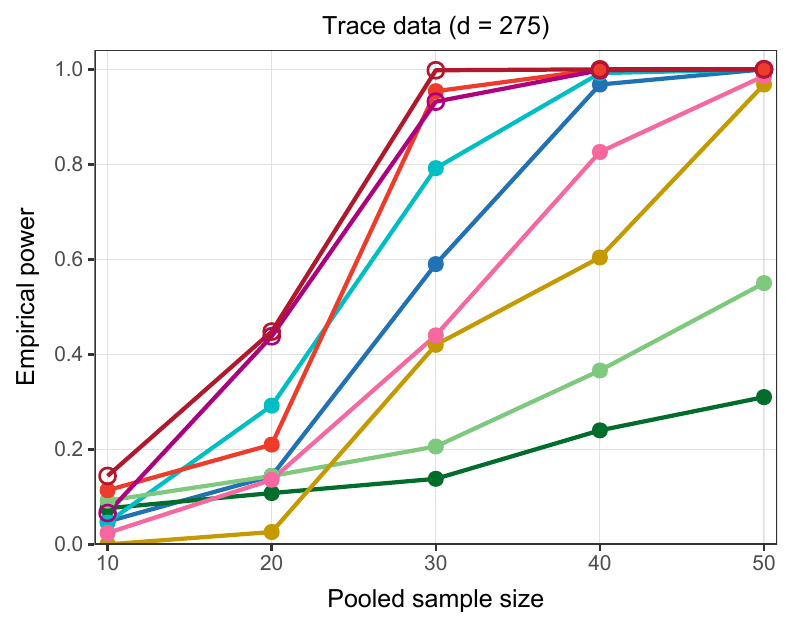}
    \end{minipage}
    \hfill
    \begin{minipage}[t]{0.31\textwidth}
        \centering
        \includegraphics[
            width=\linewidth,
            height=1.7in,
            keepaspectratio
        ]{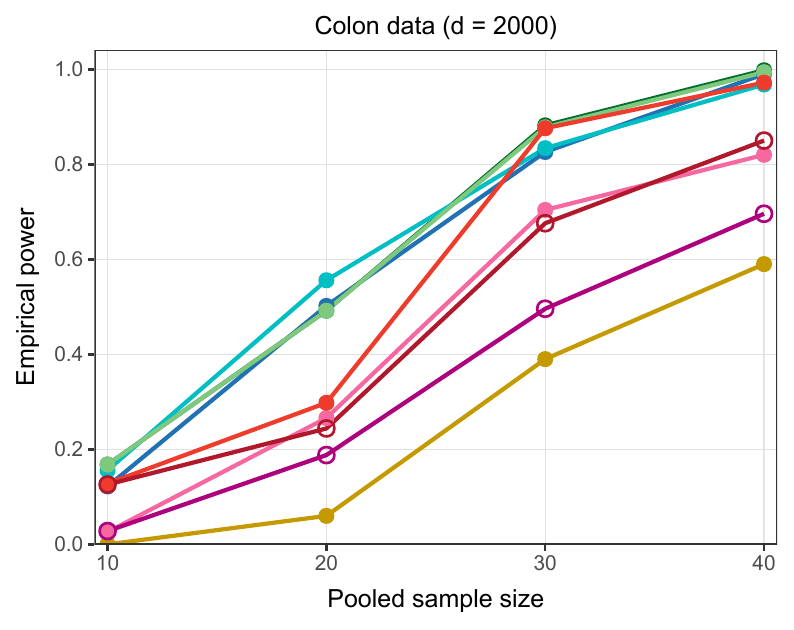}
    \end{minipage}

    \vspace{0.15cm}


    \begin{minipage}[t]{0.31\textwidth}
        \centering
        \includegraphics[
            width=\linewidth,
            height=1.7in,
            keepaspectratio
        ]{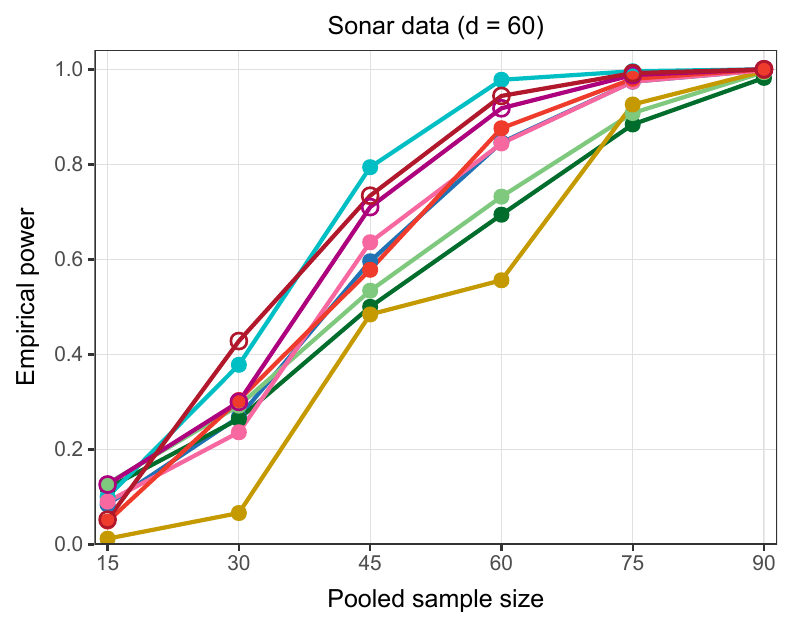}
    \end{minipage}
    \hfill
    \begin{minipage}[t]{0.31\textwidth}
        \centering
        \includegraphics[
            width=\linewidth,
            height=1.7in,
            keepaspectratio
        ]{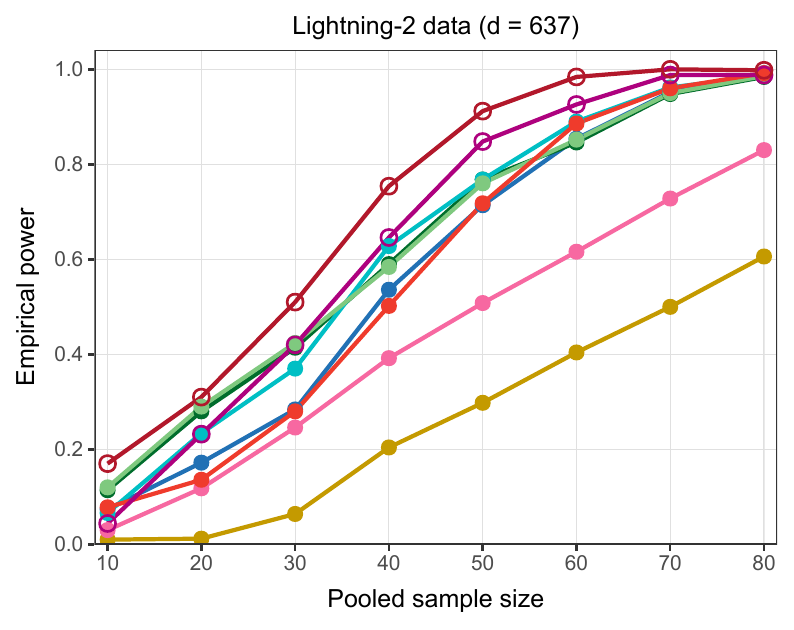}
    \end{minipage}
    \hfill
    \begin{minipage}[t]{0.31\textwidth}
        \centering
        \includegraphics[
            width=\linewidth,
            height=1.7in,
            keepaspectratio
        ]{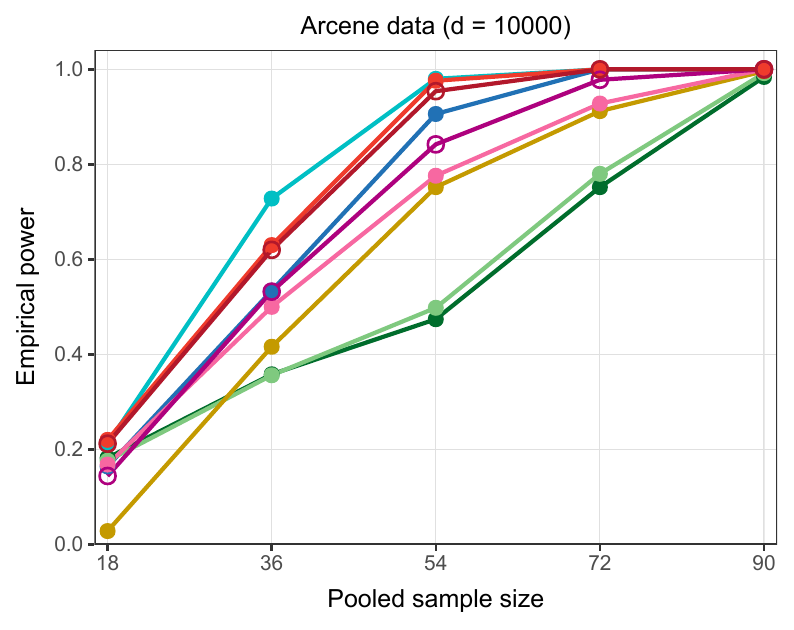}
    \end{minipage}

    \caption{
    Powers of Energy (\textcolor[HTML]{7FC97F}{$\bullet$}), MMD (\textcolor[HTML]{006D2C}{$\bullet$}), MST-run (\textcolor[HTML]{2171B5}{$\bullet$}), NN (\textcolor[HTML]{00BFC4}{$\bullet$}), MDP (\textcolor[HTML]{C49A00}{$\bullet$}), ODER (\textcolor[HTML]{F768A1}{$\bullet$}), ODER-SHC (\textcolor[HTML]{EF3B2C}{$\bullet$}), ODER-L1 (\textcolor[HTML]{AE017E}{$\circ$}) and ODER-SHC-L1 (\textcolor[HTML]{B2182B}{$\circ$}) tests in benchmark datasets.
    }
    \vspace{-0.15in}
    \label{fig:real-data-power}
\end{figure}

These datasets have been extensively used in the classification literature, and the two classes in each dataset are known to exhibit reasonable separation. Thus, we treat the alternative hypothesis as true and compare the competing tests based on their empirical powers. However, when the tests are applied to the complete data sets, most of them reject the null hypothesis of distributional equality. A comparison based on a single application to the full data is therefore not informative. To estimate the power, for each data set, we repeatedly applied the tests to 500 random subsamples and calculated the proportion of rejections of \({\cal H}_0\). We considered different subsample sizes, and the results are reported in Figure~\ref{fig:real-data-power}. For the Ionosphere, Trace, and Colon data sets (top row in Figure~\ref{fig:real-data-power}
) we used balanced subsamples, with equal numbers of observations drawn from the two classes. For the other three data sets, we used unbalanced subsamples, with the proportions of observations from the two classes approximately matching those in the original data.

Figure 9 shows the empirical powers of the different tests for these six datasets. The proposed tests generally performed well across the datasets, although their performance depends on the dataset and the choice of distance. For the Ionosphere data, the MMD and Energy tests showed relatively high power, while the proposed tests also performed well and generally outperformed their other three competitors. The MDP test had comparatively lower power, particularly for smaller pooled sample sizes. The Sonar data presented a somewhat different picture. The NN test had the best performance, followed by the proposed tests and the MST-run test. The MMD and Energy tests had somewhat lower power for larger sample sizes. As in the Ionosphere data, the MDP test had comparatively lower power.

For the Trace data, the proposed tests based on the \(\ell_1\)-distance showed high power for most sample sizes. The ODER-SHC test matched their performance for larger samples. The NN test also showed a substantial increase in power with sample size. The ODER and MST-run tests also had satisfactory performance. The MDP, Energy and MMD tests, particularly the last two, had comparatively lower power.
A somewhat similar pattern was observed for the Lightning-2 data, where the proposed tests based on the \(\ell_1\)-distance showed high power across the considered sample sizes. The ODER-SHC, NN, Energy, MST-run, and MMD tests also performed well, with power increasing steadily with sample size and approaching that of the \(\ell_1\)-distance based tests for larger sample sizes. The ODER test based on the Euclidean distance had comparatively lower power. The MDP test showed the lowest power among the methods considered here.

For the Colon data, the NN, MST-run, MMD, and Energy tests generally showed higher power than the proposed methods. Among the proposed tests, the ODER and ODER-SHC tests outperformed their corresponding \(\ell_1\)-distance based versions. Nevertheless, the proposed tests showed higher power than the MDP test.
Finally, for the Arcene data, all proposed methods performed well, with NN, ODER-SHC and ODER-SHC-L1 tests showing high power for moderate and large sample sizes. Among other competing methods, the MST-run test showed high power, followed by the ODER and ODER-L1 tests. Compared to them, the MDP test had slightly lower power. The performance of MMD and Energy tests was relatively inferior in this example.

\begin{figure}[t]
    \centering

    \begin{tabular}{cc}
        \includegraphics[height=1.6in, width=0.4\textwidth]{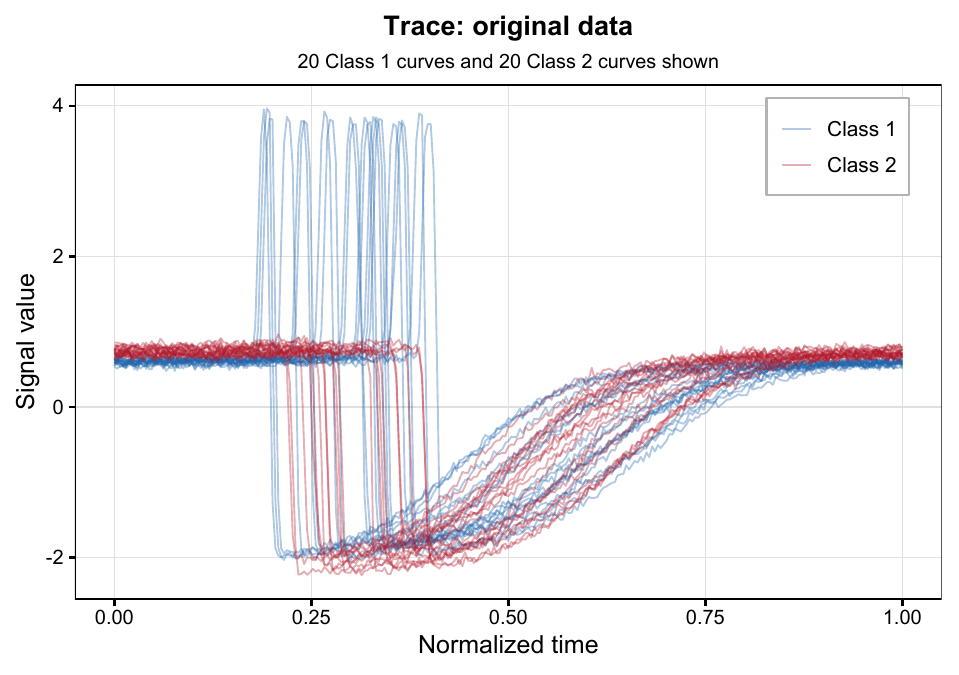} &
        \includegraphics[height=1.6in, width=0.4\textwidth]{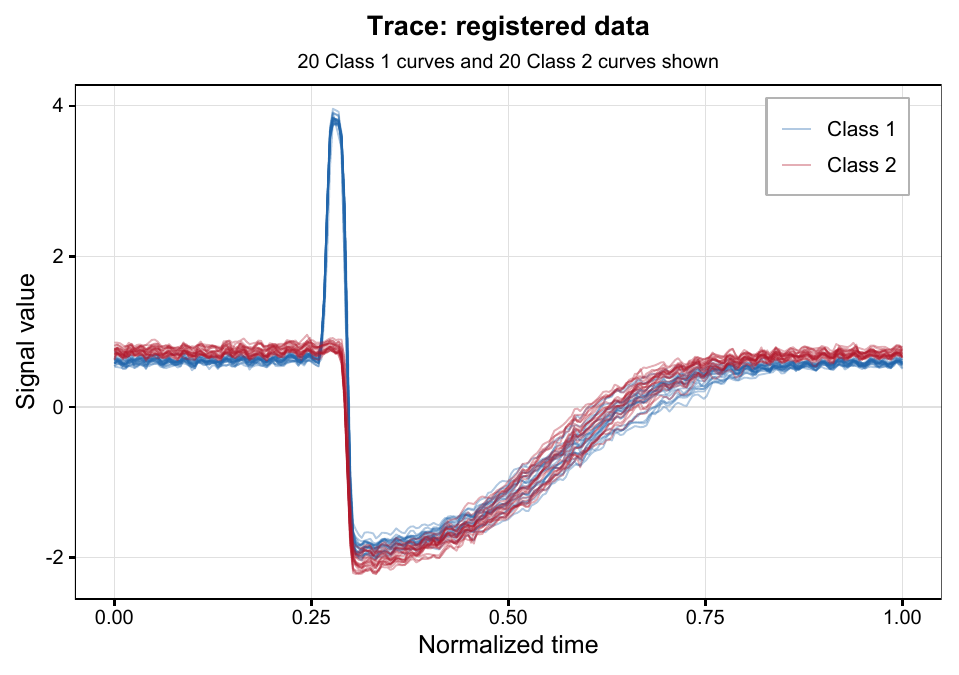} \\
     
    \includegraphics[height=1.6in, width=0.4\textwidth]{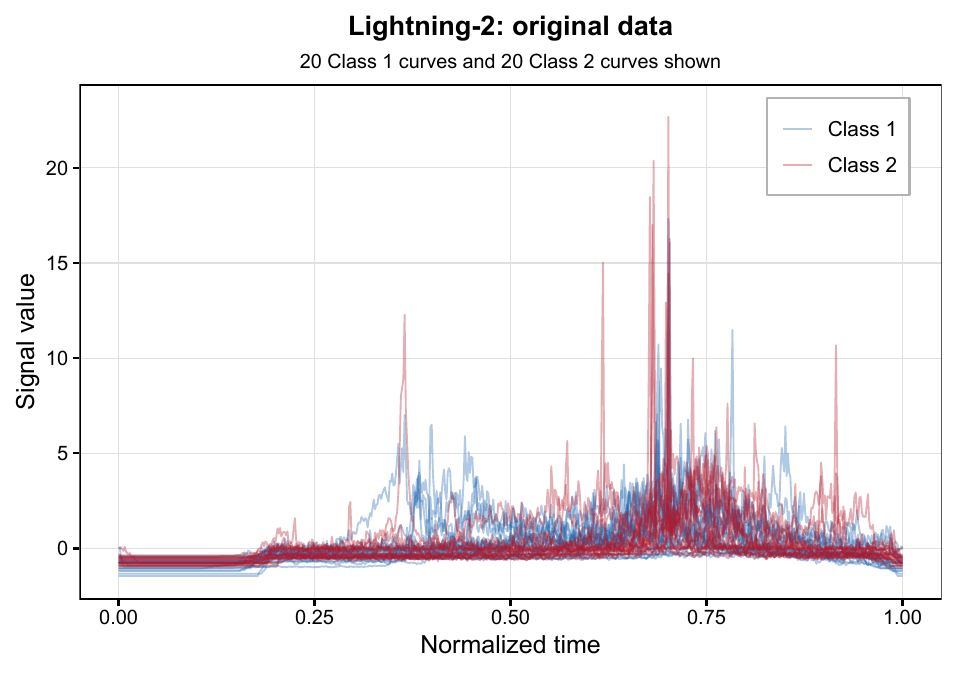} &
        \includegraphics[height=1.6in,width=0.4\textwidth]{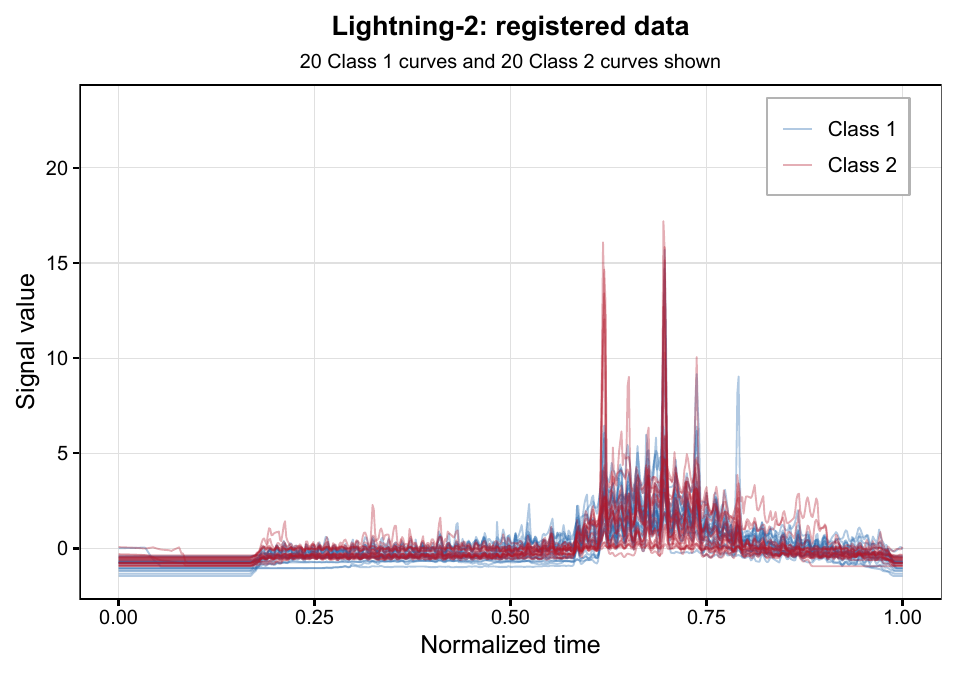}
\end{tabular}

    \vspace{0.06in}

    \caption{
    Trace data and Lightning-2 data before and after Registration
    }
    \vspace{-0.15in}
    \label{fig:registration visualization}
\end{figure}


For the time series data, dynamic time-warping distance \citep{berndt1994using} is often preferred to the Euclidean distance or $\ell_1$-distance. Another option is to go for data registration \citep{wang1997alignment,ramsay1998curve} to align the observations before using 
statistical methods based on conventional distance functions.
For Trace data and Lightning-2 data, we used the
registration method based on the square root velocity framework \citep{Tucker2013},  available in the {\tt R} package {\tt fdasrvf}. However, for the registered  Trace data, all competing tests had power $1$ even for $m=n=5$. This is explained by Figure 10, which shows clear separation between the two classes after registration. In particular, the two classes differ substantially in their locations over a particular interval, resulting in much larger interclass than intraclass distances. Consequently, all the tests performed well on the registered data. However, this was not the case for the Lightning-2 data, for which registration did not produce any visible improvement in class separability. The performances of the different methods were almost similar to those reported in Figure 9. Therefore, to avoid repetition, we do not report the results for the registered Lightning-2 data.

\section{Results from the analysis of simulated data sets}

For further evaluation of the proposed methods, we analyze some simulated datasets and compare their performance with that of the MMD, Energy, MST-run and NN tests. Although the MDP test is a special case of the derangement framework, its performance was consistently inferior in the benchmark-data analysis. In the simulated data sets also, it
always had lower power than our proposed tests. Therefore, we do not report those results here and focus instead on other competing procedures that showed more competitive performance. 
We consider eight examples with $m=n=50$ and $d=2^i$ for $i=1,2,\ldots,11$. For each setting, the experiment is repeated $500$ times, and the empirical power of each test is computed as the proportion of rejections of ${\cal H}_0$. The resulting powers are reported in Figure 11.

To examine the level properties of different tests, first we generate $100$ observations from $\mathcal{N}_d({\bf 0}_d,{\bf I}_d)$ and divide them into two equal-sized samples. We refer to this as {\bf Example 0}. Figure 11 shows that in this example, for all values of $d$, the empirical rejection rates of all tests were close to the nominal level of \(5\%\). We also repeated the experiment using distributions other than \(\mathcal{N}_d(\mathbf 0_d,\mathbf I_d)\), but our basic findings remained almost the same.

To investigate the power properties of the competing tests, next we consider two examples involving location and scale differences, but unlike Examples A and B, here we consider correlations among the coordinate variables.

{\bf Example 1:}
Consider two normal  distributions  $\mathcal{N}_d(\mathbf{0}_d,\sigmat_d)\) and \(\mathcal{N}_d(0.3\mathbf{1}_d,\sigmat_d)\) having the same dispersion matrix $\sigmat_d=((\sigma_{ij}))$ with $\sigma_{ij}=(0.5)^{|i-j|}$.  These two distributions 
differ only in their locations.

{\bf Example 2:} It also deals with
two normal distributions \(\mathcal{N}_d(\mathbf{0}_d,\sigmat_d)\) and \(\mathcal{N}_d(\mathbf{0}_d,1.5 \sigmat_d)\), where $\sigmat_d$ as in Example 1. These two distributions differ only in their scales. 

In the location problem (Example 1), the MMD and Energy tests had relatively higher power, while the powers of other competing tests were comparable. In the scale problem (Example 2), the MMD test had the best performance, closely followed by ODER-SHC, ODER-SHC-L1 and Energy tests. The ODER and ODER-L1 tests had relatively lower powers at smaller dimensions, but their powers gradually increased with $d$. 
In contrast, the NN and MST-run tests performed poorly in this example, with their powers decreasing toward zero as the dimension increased. The poor performance of these two tests in high dimensions has also been discussed by \cite{biswas2014nonparametric}.

In the next three examples, two distributions have the same mean vector and dispersion matrices with the same trace. Thus, we  have  $\nu^2 =0$ and $\sigma_1^2=\sigma_2^2$, where $\nu^2$, $\sigma_1^2$ and $\sigma_2^2$ are defined as
in Assumption (A3). In Example 6 also, we have $\nu^2 =0$ and $\sigma_1^2=\sigma_2^2$. These examples allow us to investigate the performance of the proposed and competing tests when the assumptions of Theorems 3.1 and 3.2 are not satisfied.

{\bf Example 3:} It deals with a problem involving a normal distribution $N_d({\bf 0}_d, 3{\bf I}_d)$ and the standard multivariate $t$-distribution with $3$ degrees of freedom.

{\bf Example 4:} 
Consider two normal distributions \(N_d(\mathbf{0}_d,\sigmat_\mathbf{1,d}^{\circ})\) and \({N}_d(\mathbf{0}_d,\sigmat_\mathbf{2,d}^{\circ})\) that differ only in their correlation structures. The $(i,j)$-th ($1 \le i,j \le d$) element of \(\sigmat_\mathbf{1,d}^{\circ}\) is \(0.3^{|i-j|}\), whereas that of \(\sigmat_\mathbf{2,d}^{\circ}\) is \((-0.3)^{|i-j|}\). 

{\bf Example 5:}  For $\thetavec=(\theta^{(1)},\ldots,\theta^{(d)})^\top$ and $\sigma >0$, define ${\cal L}_d(\thetavec, \sigma {\bf I}_d)$ as the distribution of a $d$-dimensional random vector $\Xvec=(X^{(1)},\ldots,X^{(d)})^{\top}$ with independent coordinates, where $X^{(i)}$ ($i=1,2,\ldots,d$) follows the univariate Laplace distribution with location $\theta^{(i)}$ and scale $\sigma$. Here, one distribution is an equal mixture of
\(\mathcal L_{d}(0.3{\mathbf 1}_d, {\bf I}_{d})\) and \(\mathcal L_{d}(-0.3\mathbf{1}_d,\, 4{\bf I}_{d})\), while the other is an equal mixture of 
\(\mathcal L_{d}(\thetavec_{2},\, {\bf I}_{d})\) and \(\mathcal L_{d}(-\thetavec_{2},\, 4{\bf I}_{d})\),
where \(\thetavec_{2}=(\theta_{2}^{(1)},\ldots,\theta_{2}^{(d)})^{\top}\) with
\(\theta_{2}^{(i)} = (-1)^{i}(0.3)\) for \(i\in\{1,\ldots,d\}\).

{\bf Example 6:} For $0\leq a<b$, let $H_d(a,b)$ be the $d$-dimensional region obtained by removing the smaller hypercube $[-a,a]^d$  from the larger hypercube $[-b,b]^d$. Here we consider two distributions, one of which is $U(H_d(1,2))$, the uniform distribution on $H_d(1,2)$, and the other one is an equal mixture of $U(H_d(0,1))$ and $U(H_d(2,\sqrt{7}))$.

\begin{figure}[!b]
    \centering

    \begin{tabular}{
        c@{\hspace{0.05in}}
        c@{\hspace{0.05in}}
        c
    }


        \includegraphics[
            width=0.29\textwidth,
            height=0.18\textheight,
            keepaspectratio
        ]
        {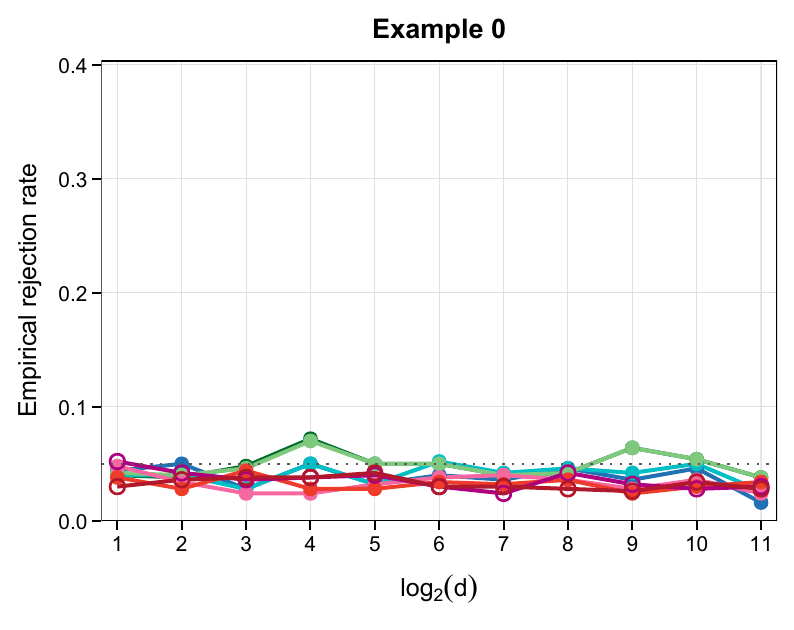}
        &
        \includegraphics[
            width=0.29\textwidth,
            height=0.18\textheight,
            keepaspectratio
        ]
        {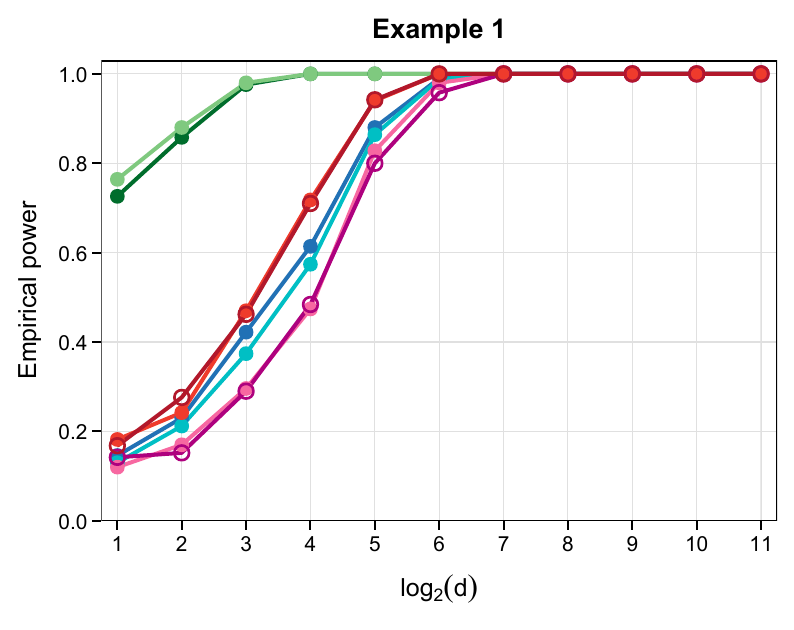}
        &
        \includegraphics[
            width=0.29\textwidth,
            height=0.18\textheight,
            keepaspectratio
        ]
        {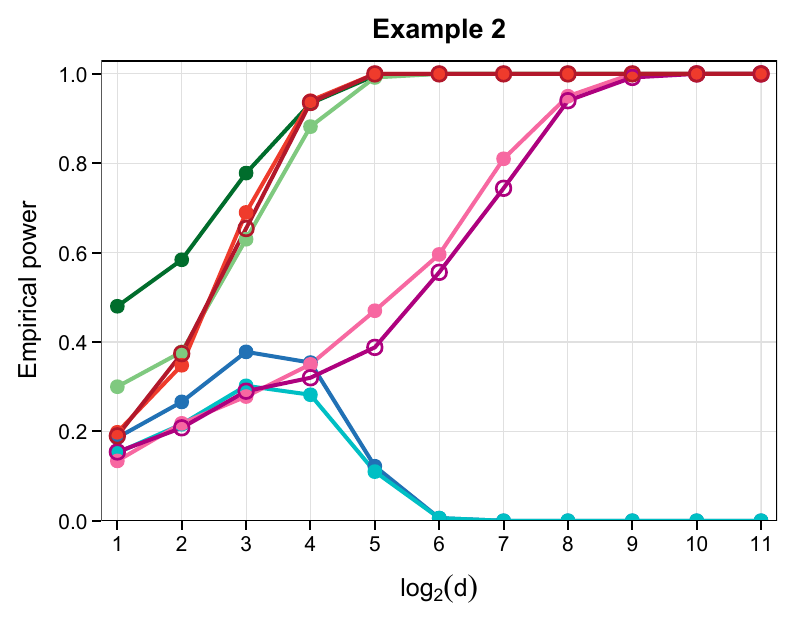}

        \\[-0.04in]


        \includegraphics[
            width=0.29\textwidth,
            height=0.18\textheight,
            keepaspectratio
        ]
        {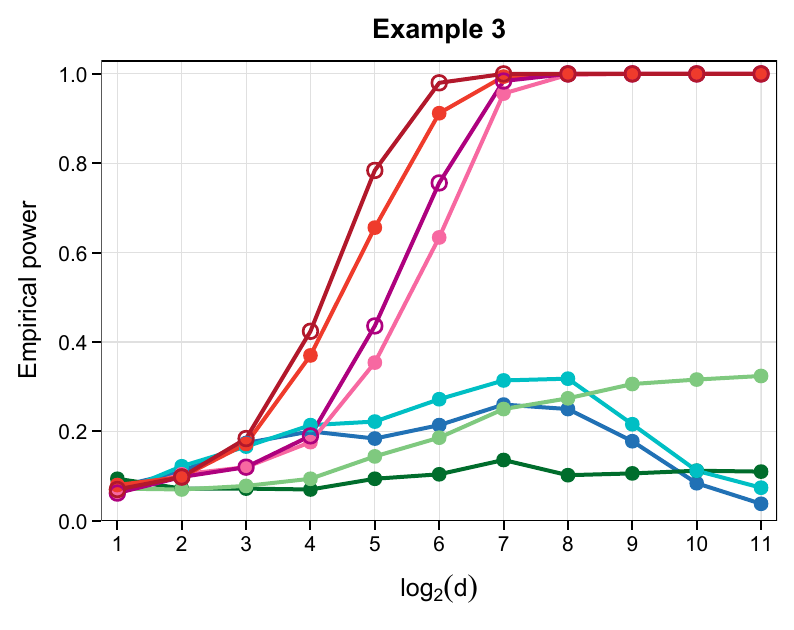}
        &
        \includegraphics[
            width=0.29\textwidth,
            height=0.18\textheight,
            keepaspectratio
        ]
        {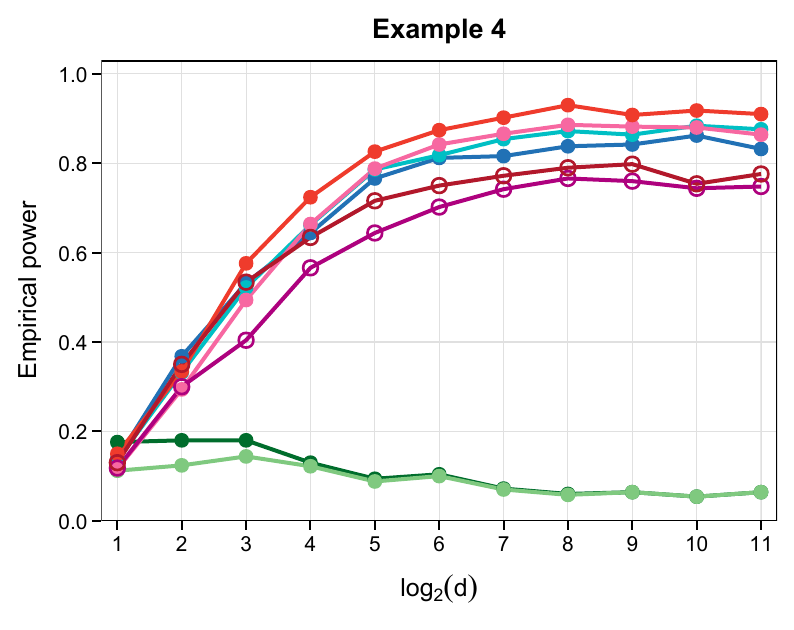}
        &
        \includegraphics[
            width=0.29\textwidth,
            height=0.18\textheight,
            keepaspectratio
        ]
        {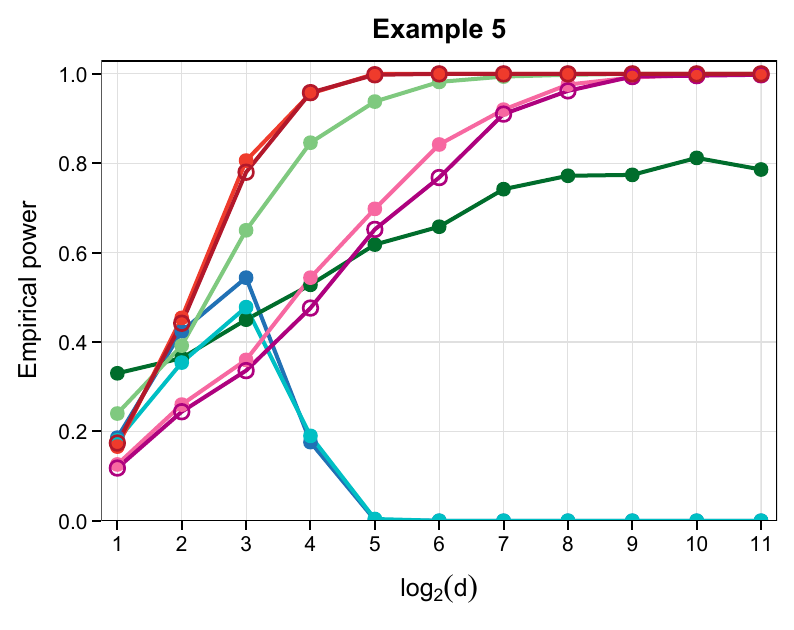}

        \\[-0.04in]


        \includegraphics[
            width=0.29\textwidth,
            height=0.18\textheight,
            keepaspectratio
        ]
        {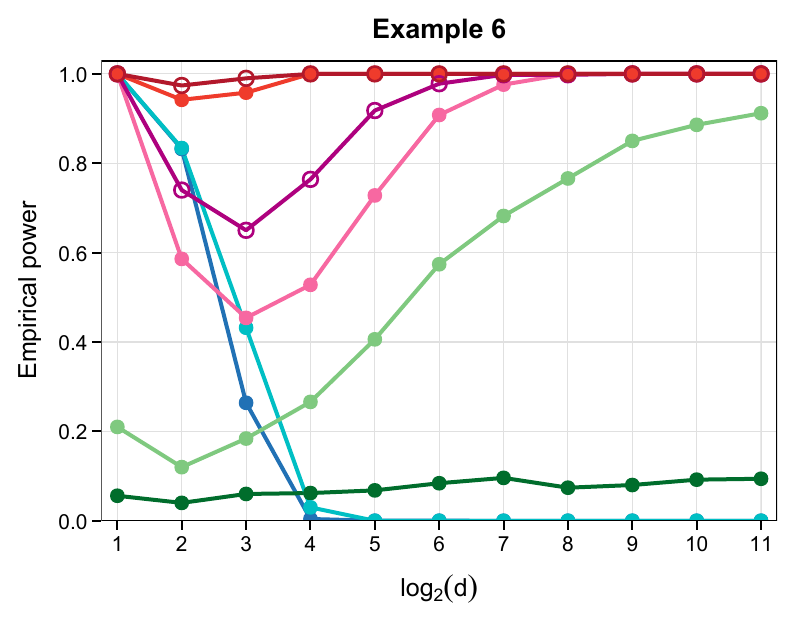}
        &
        \includegraphics[
            width=0.29\textwidth,
            height=0.18\textheight,
            keepaspectratio
        ]
        {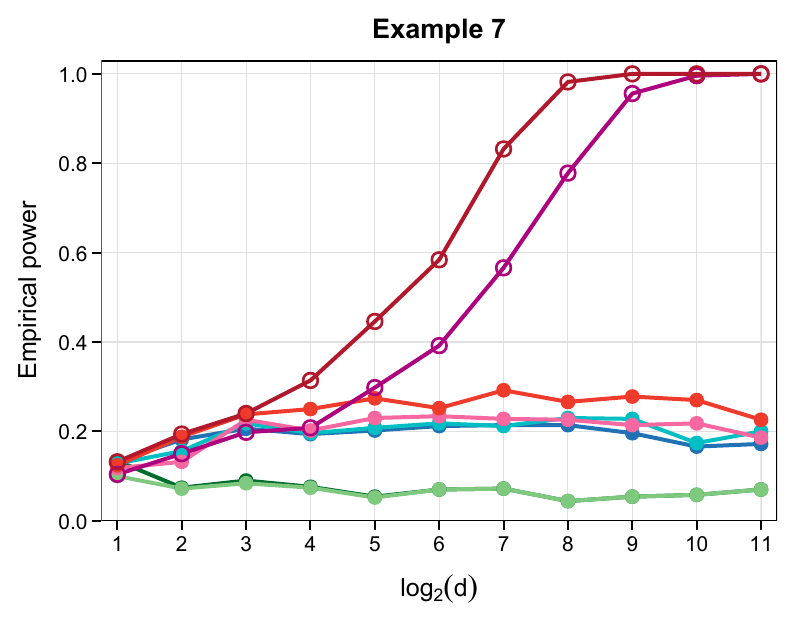}
        &
        \includegraphics[
            width=0.29\textwidth,
            height=0.18\textheight,
            keepaspectratio
        ]
        {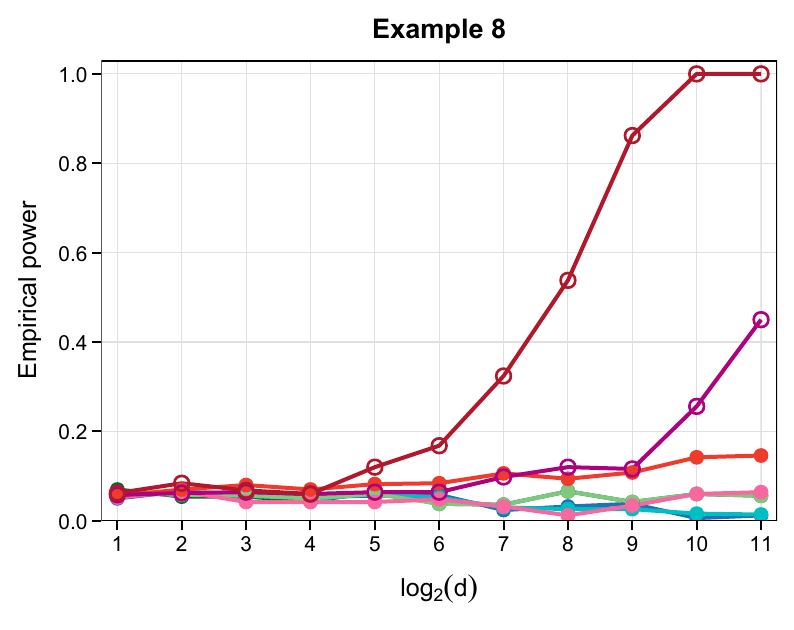}

    \end{tabular}

    \vspace{-0.03in}

    \caption{
    Empirical rejection rates in Example-0 and powers in
    Examples 1--8 for Energy
    (\textcolor[HTML]{7FC97F}{$\bullet$}),
    MMD
    (\textcolor[HTML]{006D2C}{$\bullet$}),
    MST-run
    (\textcolor[HTML]{2171B5}{$\bullet$}),
    NN
    (\textcolor[HTML]{00BFC4}{$\bullet$}),
    ODER
    (\textcolor[HTML]{F768A1}{$\bullet$}),
    ODER-SHC
    (\textcolor[HTML]{EF3B2C}{$\bullet$}),
    ODER-L1
    (\textcolor[HTML]{AE017E}{$\circ$})
    and ODER-SHC-L1
    (\textcolor[HTML]{B2182B}{$\circ$}) tests.
    }

    \label{fig:simulation-power}

\end{figure}


In Example 3, our proposed tests outperformed all other tests considered here. While the powers of all competing tests remained below 0.4, those of the proposed tests increased steadily and approached $1$ as the dimension increased. 
At least one of the conditions in each of Theorems 3.1-3.4 is violated in this example. 
Nevertheless, the proposed tests based on both the Euclidean and $\ell_1$-distances exhibited steadily increasing power with $d$.
This illustrates that the conditions of Theorems 3.1-3.4 
are sufficient but not necessary for the high-dimensional consistency of the proposed tests.
In this example, the distributions of each coordinate variable differ between the two distributions, but they have the same location and the same dispersion matrix. This may explain the slightly better performance of the \(\ell_1\)- distance based tests compared with their Euclidean distance based counterparts.
In Example 4, the MMD and Energy tests performed poorly, whereas the other tests had reasonably competitive power.  The proposed tests based on Euclidean distance had the best overall performance, closely followed by NN, MST-run and the tests based on $\ell_1$-distance.


Examples 5 and 6 deal with non-elliptic mixture distributions. In these examples, the MST-run and NN tests performed poorly, with their powers decreasing toward zero as the dimension increased. In both examples, ODER-SHC and ODER-SHC-L1 outperformed all other tests considered here. In Example 5, the Energy test had the next best performance, closely followed by the ODER and ODER-L1 tests. In Example 6, however, all the proposed tests showed higher power than their competitors for $d>4$. The MMD test performed poorly in this example.

Finally, Examples 7 and 8 are designed to resemble Examples C and D, respectively. In Examples C and D, the ODER and ODER-SHC tests showed unsatisfactory performance, whereas their \(\ell_1\)-distance based counterparts performed well. We use these examples to examine how the other competing tests perform in such situations.

{\bf Example 7:} It deals with two distributions, each having i.i.d. coordinate variables. For the first distribution, they are i.i.d.
$N(0, 5/3)$, whereas for the second, they are i.i.d. $t_5$, the standard $t$-distribution with 5 degrees of freedom.

{\bf Example 8:} Consider two multivariate normal distributions with the same mean vector ${\bf 0}_d$ but different dispersion matrices. One has dispersion matrix \(2{\bf I}_d\), while the other has a diagonal dispersion matrix whose first \(d/2\) diagonal elements are 1 and the remaining \(d/2\) elements are 3.

Figure 11 shows that the \(\ell_1\)-distance based proposed tests, particularly ODER-SHC-L1, performed well in both examples, whereas the other tests had substantially lower power. ODER and ODER-SHC tests also outperformed the competing tests, but had lower power than their \(\ell_1\)-distance based counterparts.

\section{Large sample behavior of the proposed tests}
\label{sec:large-sample}

In the previous sections, we established the consistency of the proposed tests in the HDLSS regime under appropriate regularity conditions and demonstrated their empirical performance through analyses of several benchmark and simulated data sets. We now investigate their theoretical behavior in the classical asymptotic regime, where the dimension is fixed, and the sample sizes tend to infinity. Throughout this section, we assume \(m/N\to p\in(0,1)\) as \(N=m+n\to\infty\).

First, let us study the behavior of the test statistic $D_{m,n}$
under the permutation null distribution.
Under ${\cal H}_0:F=G$, the
sample labels are exchangeable. So, conditioned on the pooled
observations $
\{\Zvec_1,\ldots,\Zvec_N\}
$
and the optimal derangement $\pi_N$, all permutations of the two types of labels ($m$ labels of one type and $n$ labels of the other type) on $\Zvec_1,\ldots,\Zvec_N$ are
equally likely. Let $L_i$ be the label attached to
the $i$-th ($i=1,2,\ldots,N$) observation and 
$
W_i
=
I\{L_i\neq L_{\pi_N(i)}\}
$ be the indicator variable corresponding to the $i$-th directed edge $\E_i=(\Zvec_i,\Zvec_{\pi_N(i)})$. 
Thus, the number of cross-sample edges is given by
$
D_{m,n}=\sum_{i=1}^N W_i.
$
Since $\pi_N(i)\neq i$, under the permutation distribution, $
\P_{\perm}
\{L_i=1,L_{\pi_N(i)}=0\}$
=
$\P_{\perm}
\{L_i=0,L_{\pi_N(i)}=1\}$
=
${mn}/{N(N-1)}.
$ 
Therefore, we have
\[
a_{m,n}:=\E_{\perm}(W_i)
=
{2mn}/{N(N-1)},
~~~~\text{  and hence  }~~~~
\nu_{m,n}:=\E_{\perm}(D_{m,n})
=
Na_{m,n}
=
{2mn}/({N-1}).
\tag{6.1}
\label{eq:null-mean}
\]

Now, we calculate the variance of $D_{m,n}$. Note that it can be expressed as 
\begin{align*}
\hspace{-0.3in}
\Var_{\perm}(D_{m,n}) 
=
\sum_{i=1}^N\Var_{\perm}(W_i)
+
2\hspace{-0.15in}\sum_{1\leq i<j\leq N}
\cov_{\perm}(W_i,W_j)
 =
N\big[a_{m,n}-a_{m,n}^2\big]
+
\sum_{i\neq j} \Big[
\E_{\perm}(W_iW_j)-a_{m,n}^2\Big].
\end{align*}
To compute $\E_{\perm}(W_iW_j)$, we consider three types of possible cases depending on the number of common vertices in  $
e_i=(\Zvec_i,\Zvec_{\pi_N(i)})$ and $
e_j=(\Zvec_j,\Zvec_{\pi_N(j)}).
$

($i$) If $\E_i$ and $\E_j$ form a cycle of size $2$, then we have
$W_i=W_j$. This leads to $\E_{\perm}(W_iW_j)
=
a_{m,n}.
$

($ii$)
If $\E_i$ and $\E_j$ have exactly one common vertex, $W_iW_j$ takes the value $1$ if and only if the label of the common vertex differs from the labels of the other two vertices. So, in this case
\[
\begin{aligned}
\E_{\perm}(W_iW_j)
=
\frac{
m(m-1)n+n(n-1)m
}{
N(N-1)(N-2)
}
=
\frac{mn}{N(N-1)}
=\frac{1}{2}a_{m,n}.
\end{aligned}
\]
($iii$) If $\E_i$ and $\E_j$ have no common vertex, it is easy to see that
\[
\E_{\perm}(W_iW_j)
=
\frac{
4m(m-1)n(n-1)
}{
N(N-1)(N-2)(N-3)
}
=:b_{m,n}.
\]

Let $C_N$ be the number of cycles of size 2 in the optimal derangement graph. Note that every vertex belonging to a cycle of length at least three gives
a pair of edges (the incoming and outgoing edges at that vertex) of the second type. So, in such cases, the numbers of pairs of the first and second types are $C_N$ and $N-2C_N$, respectively. Therefore,
the number of unordered pairs of the third type turns out to be
$\binom{N}{2}-N+C_N$.
Combining these three cases, we get

\vspace{-0.2in}
\[
\begin{aligned}
\Var_{\perm}(D_{m,n})
={}
N\bigl(a_{m,n}-a_{m,n}^2\bigr)
&+
2C_N
\bigl(a_{m,n}-a_{m,n}^2\bigr)
+
2\bigl(N-2C_N\bigr)
\bigl(\frac{1}{2}a_{m,n}-a_{m,n}^2\bigr)\\
&+
2\left\{
\binom{N}{2}-N+C_N
\right\}
\bigl(b_{m,n}-a_{m,n}^2\bigr)
\end{aligned}
\tag{6.2}
\vspace{-0.1in}
\label{eq:null-variance}
\]
This can be further simplified, and for $N\geq4$, we have
\[
\omega_{m,n}^2:=
\Var_{\perm}(D_{m,n})=
\frac{
4mn(m-1)(n-1)
\left\{
N(N-3)+2C_N(N-1)
\right\}
}{
N(N-1)^2(N-2)(N-3)
}.
\tag{6.3}
\label{eq:null-variance-simplified}
\]

Note that while \(\nu_{m,n}\) depends only on sample sizes, 
\(\omega_{m,n}^2\) depends on 
\(\pi_N\), or equivalently, on the topology of the optimal derangement graph. This is analogous to the MST run test of \cite{MST}, whose conditional null variance depends on the topology of the observed MST. However, once the pooled observations and \(\pi_N\) are given, \(C_N\) is fixed and directly observable.
If the optimal derangement graph contains cycles of size 2 only, then $C_N= N/2$, and $\omega_{m,n}$ gives the null variance for the MDP test of \cite{MinDist}. If it consists of a single cycle of size $N$, then $C_N=0$, and $\omega_{m,n}$ gives the variance of $D_{m,n}^{\circ}$. So, \(D^\circ_{m,n}\) usually has smaller variance than \(D_{m,n}\). This smaller stochastic variation may partly explain the better performance of ODER-SHC than ODER in most of our examples. Conversely, the higher variance of the MDP statistic may partly explain its relatively lower performance observed in Section 4.

As discussed earlier, the ODER-SHC test is distribution-free, whereas the ODER test requires the permutation method to determine the critical value. This can be computationally expensive for large samples. In such cases, it is preferable to find the critical value based on the asymptotic null distribution. Given the pooled observations and the optimal derangement, we have already derived the mean \(\nu_{m,n}\) and variance \(\omega_{m,n}^2\) of the permutation null distribution. The large-sample permutation distribution of \(D_{m,n}\) is given by the following theorem. 


\begin{thm}
\label{thm:null-asymptotic}
If $m/N\rightarrow p \in (0,1)$ as $N \rightarrow \infty$, then 
conditional on the pooled observations and the optimal
derangement,
\(
({D_{m,n}-\nu_{m,n}})/{\omega_{m,n}}
\xrightarrow{d}
N(0,1),
\)
where $\nu_{m,n}$ and $\omega_{m,n}^2$ are as defined in
\eqref{eq:null-mean} and \eqref{eq:null-variance-simplified}, respectively.
\end{thm}

The asymptotic permutation distribution of $D_{m,n}$ derived in Theorem 6.1 applies under both ${\cal H}_0$ and ${\cal H}_1$ since its derivation depends only on the combinatorial structure of the derangement and the permutation of the sample labels, not on the equality of the underlying distributions (see the proof in the Appendix).
The asymptotic null distribution of $D_{m,n}^{\circ}$ is obtained from it by setting $C_N=0$ in the expression of $\omega_{m,n}$. These results also hold for the corresponding tests based on $\ell_1$-distances (see Remark 5 in the Appendix). 

Table 1 reports the ratio of the computing time required by the permutation method to that required to obtain the critical value from the large-sample distribution. The computational advantage of the latter is evident from the table.

Table 2 shows the empirical rejection rates of the ODER tests, based on the (i) permutation method and (ii) large-sample distribution, for different choices of $d$ and $N$, calculated from 5000 repetitions. When both distributions are normal or both are Laplace, for all choices of $d$ and $N$, the large-sample method rejected \({\cal H}_0\) in almost \(5\%\) of the cases, while the permutation method was slightly conservative.  But this mild conservativeness did not have any visible effect on its power when we discriminate between the standard normal and standard Laplace distributions. In almost all cases, the permutation method had nearly the same power as the large sample method. 

Next comes the question of large-sample consistency. For that study, we assume that

\begin{enumerate}
    
\item[(A4)] $F$ and $G$ are absolutely continuous with respect to the $d$-dimensional Lebesgue measure, and they have densities $f$ and $g$, respectively. 

\item[(A5)] $\E_F\|\Xvec\|^r+\E_G\|\Yvec\|^r
<\infty$ 
for some
$
r>\frac{d}{d-1}
$.
\end{enumerate}

\begin{table}[!b]
    \caption{Ratio of the time taken by the permutation method to that of the large sample method for the computation of the critical value for the ODER test.}
    \label{tab:time}
    \centering
    \begin{tabular}{ccccccc} \hline
       N  & 100 & 200 & 400 & 1000 & 2000 \\
      Time ratio   & 105.32& 222.22 & 438.61 & 1333.33 & 2285.71 \\ \hline
    \end{tabular}
\end{table}


\begin{table}[t]
    \caption{Empirical rejection rates of the ODER tests based on permutation methods and large-sample distribution.}
    \label{tab:asymptotic-size}
    \centering
    \setlength{\tabcolsep}{4pt}
    \begin{tabular}{cc@{\hspace{0.5cm}}ccccc@{\hspace{0.5cm}}ccccc}
        \hline
        & & \multicolumn{5}{c}{$d=2$}
        & \multicolumn{5}{c}{$d=20$} \\
        \cline{3-12}
        & $N$
        & 100 & 200 & 400 & 1000 & 2000
        & 100 & 200 & 400 & 1000 & 2000 \\
        \hline

        $F$: $N_d({\bf 0}_d,{\bf I}_d)$
        & Permutation
        & 0.0360 & 0.0422 & 0.0438 & 0.0466 & 0.0500
        & 0.0400 & 0.0346 & 0.0398 & 0.0488 & 0.0516 \\

        $G$: $N_d({\bf 0}_d,{\bf I}_d)$
        & Large-sample
        & 0.0404 & 0.0504 & 0.0556 & 0.0502 & 0.0504
        & 0.0504 & 0.0538 & 0.0434 & 0.0542 & 0.0566 \\
        \hline

        $F$: ${\cal L}_d({\bf 0}_d,{\bf I}_d)$
        & Permutation
        & 0.0370 & 0.0386 & 0.0368 & 0.0420 & 0.0504
        & 0.0372 & 0.0368 & 0.0414 & 0.0444 & 0.0462 \\

        $G$: ${\cal L}_d({\bf 0}_d,{\bf I}_d)$
        & Large-sample
        & 0.0434 & 0.0494 & 0.0470 & 0.0474 & 0.0534
        & 0.0470 & 0.0524 & 0.0452 & 0.0492 & 0.0490 \\
        \hline

        $F$: $N_d({\bf 0}_d,{\bf I}_d)$
        & Permutation
        & 0.0738 & 0.1070 & 0.1694 & 0.3368 & 0.5472
        & 0.0626 & 0.0910 & 0.1742 & 0.5446 & 0.9230 \\

        $G$: ${\cal L}_d({\bf 0}_d,{\bf I}_d)$
        & Large-sample
        & 0.0812 & 0.1296 & 0.2022 & 0.3510 & 0.5556
        & 0.0762 & 0.1228 & 0.1838 & 0.5622 & 0.9298 \\
        \hline
    \end{tabular}
\end{table}

Assumption (A5) is a mild tail condition satisfied by a broad class of distributions. 
For studying the large-sample behavior of the proposed tests, we use the general result (Theorem 1) of \cite{Arias-Castro2016}, whose proof follows from \cite{HenzePenrose1999}. It applies to the usual two-sample setup with $m/N \rightarrow p \in (0,1)$ and to any graph constructed on the pooled sample that satisfies three conditions: bounded in- and out-degrees of all vertices, convergence of the out-degree of a vertex to a positive constant, and the short-edge property (i.e., long edges are rare). The last two conditions are verified for an auxiliary model \citep{HenzePenrose1999}, where for any fixed $N$, the pooled observations $\{\Zvec_{1:N},\Zvec_{2:N}\ldots,\Zvec_{N:N}\}$ are assumed to be i.i.d. with the density $\phi_N=(mf+ng)/N$, and  conditional on an observation $\zvec$, it gets label $\Xvec$ with probability $mf(\zvec)/\{mf(\zvec)+ng(\zvec)\}$ and
label $\Yvec$ with probability $ng(\zvec)/\{mf(\zvec)+ng(\zvec)\}$.
Under this model, \cite{Arias-Castro2016} proved the short-edge property of minimum distance pairings and the large-sample consistency of the MDP test \citep{MinDist}.
Following their work, we can show that (see the proof of Theorem~\ref{thm:large-sample-consistency}) both under ${\cal H}_0: F =G$ and ${\cal H}_1:F \neq G$,
\[
\frac{D_{m,n}}{N} \stackrel{a.s.}{\rightarrow} \Lambda_p(F,G)
=
2p(1-p)
\int_{\mathbb R^d}
\frac{f(z)g(z)}
{pf(z)+(1-p)g(z)}
\,dz \text{   as  } N \rightarrow \infty.
\]
 Under ${\cal H}_0$, this limiting value turns out to be $2p(1-p)$. Now one can  check that
\[
2p(1-p)-\Lambda_p(F,G)
=
2p^2(1-p)^2
\int_{\mathbb R^d}
\frac{\{f(z)-g(z)\}^2}
     {pf(z)+(1-p)g(z)}\,dz ~\ge 0,
\]
with equality if and only if $F=G$. Thus, the limiting value of $D_{m,n}/N$ under ${\cal H}_1$ is strictly 
smaller than that under ${\cal H}_0$. This proves the large-sample consistency of the resulting test. The same result holds for the test based on $D_{m,n}^{\circ}$.

\begin{thm}
\label{thm:large-sample-consistency}

Suppose that $F$ and $G$ satisfy assumptions (A4)-(A5), and $m/N \rightarrow p \in(0,1)$ as $N \rightarrow \infty$. Then, for any fixed $d\geq2$, the powers of ODER and ODER-SHC tests converge to $1$ as $N$ tends to infinity.
\end{thm}

Although Theorem 6.2 is stated for tests based on Euclidean distances, this result holds for the corresponding tests based on the $\ell_1$-distance as well (see Remark 6 in the Appendix). The large sample consistency of the ODER-SHC test continues to hold even when the SHC is approximated by the cheapest-link algorithm. This result is stated as Lemma A.7 in the Appendix. Because all norm-induced distances are equivalent in finite dimensions, the same consistency result can be proved for the implemented version of the ODER-SHC-L1 test as well.


\section{Concluding remarks}

This article investigates high-dimensional two-sample testing from the perspective of optimal transport. 
We showed that the high-dimensional behavior of the transport-based tests can depend critically on both the reference distribution and the transportation cost. In particular, the conventional combination of a uniform reference distribution and squared Euclidean distance cost may lead to a loss of discriminatory power due to concentration of transportation costs, whereas suitable choices of the reference distribution or transportation cost can substantially improve performance. 

Motivated by this observation, we introduced a reference-free framework based on optimal derangements. The framework also provides a connection between optimal transport and existing graph-based test procedures like the MDP test \citep{MinDist}. We proposed a test based on unrestricted derangements (ODER) and its distribution-free variant (ODER-SHC) obtained by restricting the derangements to a single Hamiltonian cycle. Under appropriate assumptions, we established the consistency of these proposed tests in both high-dimensional and large-sample settings, and also for moderately  sparse alternatives. We also considered corresponding
tests based on $\ell_1$-distances, which can detect differences in marginal distributions beyond differences in the first two
moments. Our analyses of simulated and benchmark data sets show that the proposed methods are competitive with
existing two-sample tests across a variety of alternatives and data settings. Although no single procedure uniformly
dominates in all situations, the derangement-based tests often perform well in high dimensions.  

Several directions for future research remain open. First, it would be natural to extend the derangement framework to \(k\)-sample problems and develop corresponding distribution-free procedures. Although straightforward extensions of the ODER and ODER-SHC tests are possible, their theoretical properties and empirical performance warrant further investigation. Second, a more complete theory for regimes in which both the dimension and the sample size increase simultaneously would provide a bridge between the HDLSS and classical asymptotic settings considered here. Third, we have considered the Euclidean and \(\ell_1\) distances as transportation costs. It would be interesting to investigate other cost functions and determine whether they can provide greater robustness or power while retaining the distribution-free property and the large-sample and HDLSS consistency of the resulting tests. Finally, developing adaptive procedures that select or combine dissimilarity measures according to the characteristics of the data and the alternatives of interest would be of considerable practical interest.

\appendix
\section{Proofs and mathematical details}

\subsection{Results on high-dimensional behavior of the derangement-based tests}




\begin{proof}[{\bf Proof of Theorem~\ref{thm:HDLSS-derangement}}]

Under Assumptions~\textnormal{(A1)--(A3)}, as $d\rightarrow \infty$,
${\|\Xvec_i-\Xvec_j\|}/{\sqrt d}
\stackrel{P}{\longrightarrow}\sigma_1\sqrt{2}$ for $1\le i<j \le m$, 
${\|\Yvec_i-\Yvec_j\|}/{\sqrt d}
\stackrel{P}{\longrightarrow}\sigma_2\sqrt{2}$ for $1\le i<j \le n$
and
${\|\Xvec_i-\Yvec_j\|}/{\sqrt d}
\stackrel{P}{\longrightarrow}
\sqrt{\sigma_1^2+\sigma_2^2+\nu^2}$ for $1\le i\le m,~ 1\le j\le n$.
Let these three limiting values be denoted by $a$, $b$ and $c$, respectively. Note that
$\Delta=2c-a-b\ge 0$, where the 
equality holds if and only if
$\sigma_1^2=\sigma_2^2$
and
$\nu^2=0$.
For a derangement $\pi$, let
$T(\pi)$ be the number of directed cross-sample edges and 
$\mathcal C_d(\pi)
=
\sum_{i=1}^{m+n}\|Z_i-Z_{\pi(i)}\|$
be the total cost.
Since $m$ , $n\geq 2$ are fixed, 
$${\mathcal C_d(\pi)}/{\sqrt d}
=
\Big(m-\frac{1}{2}{T(\pi)}\Big)a
+
\Big(n-\frac{1}{2}{T(\pi)}\Big)b
+
T(\pi)c
+
o_P(1)=
ma+nb
+
\frac{1}{2}{T(\pi)}\Delta
+
o_P(1)
$$
Since the number of admissible derangements is finite, the above \(o_P(1)\) term can be taken uniformly over \(\pi\).\\
(a) Since
$ma+nb$
does not depend on $\pi$ and $\Delta>0$, with probability tending to one, the cost will be minimum when $T(\pi)$ takes its smallest possible value $0$. 
 Thus,
$\P(D_{m,n}=0)\longrightarrow 1$
as $d \rightarrow \infty$. Since $D_{m,n}$ only takes integer values, we have $\P(D_{m,n}=0) \rightarrow 1$
as $d$ tends to infinity.

Now we compute the probability of this event under the permutation distribution. Suppose that the selected optimal derangement has $K$ cycles of size
$\ell_1,\ldots,\ell_K$. Since every cycle is of size at least $2$, 
$
K \leq \left\lfloor {N}/{2} \right\rfloor$. Under permutation, the $r=\min\{m,n\}$ labels of the smaller sample are
uniformly assigned to $r$ of the $N$ vertices. The test statistic takes the value $0$ if and only if every
cycle is monochromatic under permutation. Hence, the $r$ vertices of the smaller
sample must be a union of some complete cycles. If this union contains $j$ cycles, then $r\geq 2j$, so
$j\leq \left\lfloor r/2\right\rfloor$. Therefore, the number of
favorable labelings is at most
\[
\sum_{j=1}^{\lfloor r/2\rfloor}\binom{K}{j}
\leq
\sum_{j=1}^{\lfloor r/2\rfloor}
\binom{\lfloor N/2\rfloor}{j}=\binom{N}{r} B_{m,n}.
\]
Since there are $\binom{N}{r}$ equally likely labelings under permutation,
\(
\mathbb{P}_{\mathrm{perm}}(D_{m,n}=0)
\leq B_{m,n}.
\)
So, if $B_{m,n}<\alpha$, the critical value for the permutation test must be bigger than $0$, and that leads to the rejection of ${\cal H}_0$. Therefore, under the given condition, the power of the test converges to $1$ as $d$
tends to infinity.

(b) For ODER-SHC, the derangement is restricted to a single
cycle of size $m+n$. Any such cycle contains at least two directed cross-sample edges. Hence, the cost is minimum when $T(\pi)=2$, and this implies
$\P(D_{m,n}^{\circ}=2)\rightarrow1$ or $\P(D_{m,n}^{\circ}=2) \rightarrow 1$ as $d \rightarrow \infty$.
Note that $D_{m,n}^{\circ}=2$ if and only if all observations from each class appear consecutively. Under ${\cal H}_0$, when all permutations on the cycle are equally likely, this happens with probability ${m!n!}/{(m+n-1)!}.$
So, if
${m!n!}/{(m+n-1)!}<\alpha$,
the ODER-SHC test rejects ${\cal H}_0$ with probability converging to one as $d$ tends to infinity.
\end{proof}

\begin{rem}
{\text   Define  $M_{11}=\max_{1\le i<j \le m}\|\Xvec_i-\Xvec_j\|$, $M_{22}=\max_{1\le i<j \le n}\|\Yvec_i-\Yvec_j\|$ and $M_{12}=\min_{i,j}\|\Xvec_i-\Yvec_j\|$. When $m$ and $n$ are finite,  under (A1)-(A3), one can show that $M_{11}/\sqrt{d}$, $M_{22}/\sqrt{d}$ and $M_{12}/\sqrt{d}$ converge in probability to $a=\sigma_1\sqrt{2}$, and $b=\sigma_2\sqrt{2}$ and $c=\sqrt{\sigma_1^2+\sigma_2^2+\nu^2}$, respectively. Under the conditions of Theorem 3.1,  $2c-a-b>0$. Consequently, $min\{a,b\}<c$ and hence $\P\big((M_{11}<M_{12})\cup(M_{22}<M_{12})\big) \rightarrow 1$
as $d \rightarrow \infty$. Thus, either all $\Xvec\Xvec$-type edges or all $\Yvec\Yvec$-type edges turn out to be smaller than all $\Xvec\Yvec$-type edges with probability tending to $1$. Without loss of generality, let us assume that all $\Xvec\Xvec$-edges are smaller than $\Xvec\Yvec$-edges. So, before considering the first $\Xvec\Yvec$-edge, the cheapest-link algorithm processes all $\Xvec\Xvec$-edges to form a path containing all $\Xvec$’s. Because of the degree constraint, only two end-points of that path can then accept cross-sample edges. Therefore, the resulting Hamiltonian cycle has exactly two cross-sample edges.
Consequently, the cheapest-link implementation of the ODER-SHC test has the high-dimensional consistency under the conditions of Theorem 3.1.
}
\end{rem}

\begin{lem}
\label{lem:growing-sample}
Suppose that $\Xvec_1,\Xvec_2,\ldots,\Xvec_m$ are independent copies of $\Xvec=(X^{(1)},\ldots,X^{(d)})^\top\sim F$ and $\Yvec_1,\Yvec_2,\ldots,\Yvec_n$ are independent copies of $\Yvec=(Y^{(1)},\ldots,Y^{(d)})^\top\sim G$, where $F$ and $G$ satisfy the conditions of Theorem 3.1.
If the centered coordinate variables
$X^{(q)}-\mathbb{E}X^{(q)}$ and
$Y^{(q)}-\mathbb{E}Y^{(q)}$ $(q=1,2,\ldots)$ are uniformly sub-Gaussian and
$m_0$-dependent (under some permutation), then the consistency results in parts (a) and (b)
of Theorem 3.1 continue to hold when $N=m+n$ increases with $d$,
provided
$({\log N})/{d}\longrightarrow 0$
as $N,d\to\infty$.\end{lem}

\begin{proof}
Consider any two independent observations ${\bf U}$ and ${\bf V}$ from the pooled sample of size $N=m+n$. For each coordinate $q=1,2,\ldots,d$, write 
\[U^{(q)}-V^{(q)}=\delta_q + \xi_q, \text{  where   } \delta_q=E(U^{(q)})-E(V^{(q)}) \text {  and    } \xi_q=\left[U^{(q)}-V^{(q)}\right]-\left[E(U^{(q)})-E(V^{(q)})\right].\]
Therefore, $\|{\bf U}-{\bf V}\|^2=\sum_{q=1}^{d} (\delta_q+\xi_q)^2$,
and $\|{\bf U}-{\bf V}\|^2-E\|{\bf U}-{\bf V}\|^2=\sum_{q=1}^d\big(\xi_q^2-E(\xi_q^2)\big)+2\sum_{q=1}^{d}\delta_q\xi_q$.
Since the centered coordinate differences are uniformly sub-Gaussian and $m_0$-dependent under some permutation, so are the random variables $\xi_q$. Therefore, $\xi_q^2 - E(\xi_q^2)$ are uniformly sub-exponential and $m_0$-dependent under some permutation.
So, for any fixed $\epsilon>0$, a Bernstein inequality for \(m_0\)-dependent sub-exponential variables gives
\(
P\Big(
\big|
d^{-1}\sum_{q=1}^{d}
\big(\xi_q^2-E\xi_q^2\big)
\big|
>{\epsilon}/{2}
\Big)
\leq C_1 e^{-c_1 d}
\)
for some \(C_1,c_1>0\), uniformly over all choices of \({\bf U}\) and \({\bf V}\).

For a within-sample pair, $\delta_q=0$ for all $q=1,2,\ldots,d$. Therefore, the above exponential bound holds for $d^{-1}\|{\bf U}-{\bf V}\|^2-d^{-1}E\|{\bf U}-{\bf V}\|^2$, uniformly over all such pairs.
For a cross-sample pair, from Assumption (A3), we have
$
\sum_{q=1}^d \delta_q^2=O(d)
$. Since the \(\xi_q\)'s are uniformly sub-Gaussian and \(m_0\)-dependent under some permutation, the weighted sum
$ \sum_{q=1}^d\delta_q\xi_q $
is also sub-Gaussian (note that under sub-Gaussianity and $m_0$ dependence of the  \(\xi_q\)'s, $\|\xi_q\|_{\psi_2} \le K$ uniformly in $q,d$ implies $\|\sum_{q=1}^d\delta_q\xi_q\|_{\psi_2} \le CK (\sum_{q=1}^{d} \delta_q^2)^{1/2}$ for some $C>0$). Hence,
$\P\Big(
\big|
\frac{2}{d}\sum_{q=1}^d\delta_q\xi_q
\big|>\frac{\epsilon}{2}
\Big)
\leq C_2e^{-c_2d}$
for some constants $C_2, c_2>0$, again uniformly over all such pairs. Combining these two exponential bounds, 
for every fixed $\epsilon>0$, we get
\(
P\Big(
\big|
d^{-1}\|{\bf U}-{\bf V}\|^2-d^{-1}E\|{\bf U}-{\bf V}\|^2
\big|>\epsilon
\Big)
\leq C\exp(-c_\epsilon d),
\)
for some constants $C,c_\epsilon>0$, uniformly over all within-sample and cross-sample pairs. Since there are at most $N^2$ pairs, a union bound gives 
$$\mathbb P\Big(
\max_{i,j}
\left|
d^{-1}\|\mathbf \Zvec_i-\mathbf \Zvec_j\|^2
-d^{-1}\mathbb E\|\mathbf \Zvec_i-\mathbf \Zvec_j\|^2
\right|>\epsilon
\Big)
\le C N^2\exp(-c_\epsilon d).
\vspace{-0.1in}
$$
Thus, $ \max_{i,j}
\left|
d^{-1}\|\mathbf \Zvec_i-\mathbf \Zvec_j\|^2
-d^{-1}\mathbb E\|\mathbf \Zvec_i-\mathbf \Zvec_j\|^2
\right| \stackrel{P}{\rightarrow} 0$
when
$N^2=o\left(\exp(c_\epsilon d)\right)$, which is satisfied if
$
({\log N})/{d}\rightarrow 0.
$
 So,
${\|\Xvec_i-\Xvec_j\|^2}/{d}
\stackrel{P}{\rightarrow} a^2$,
${\|\Yvec_i-\Yvec_j\|^2}/{d}
\stackrel{P}{\rightarrow}b^2$ 
and
${\|\Xvec_i-\Yvec_j\|^2}/{d}
\stackrel{P}{\rightarrow}
c^2$ uniformly over all pairs, 
where $a=\sigma_1\sqrt{2}$, 
$b=\sigma_2\sqrt{2}$ and $c=\sqrt{\sigma_1^2+\sigma_2^2+\nu^2}$ as defined in the proof of Theorem 3.1. Clearly, we have $2c-a-b>0$. 

Suppose that the derangement graph has two cross-sample edges $\Xvec_i\Yvec_k$ and $\Xvec_j\Yvec_{\ell}$. Now,  
($i$) if $i\neq j$ and $k \neq \ell$, replacing
them by $\Xvec_i\Xvec_j$
and $\Yvec_k\Yvec_{\ell}$ reduces the cost of the derangement. Again, ($ii$) if $i=j$ and $k \neq \ell$, choose an edge $\Xvec_s\Xvec_t$ where $s\neq t \neq i$. Replacing $\Xvec_i\Yvec_k$,
$\Xvec_i\Yvec_\ell$ and $\Xvec_s\Xvec_t$ by $\Xvec_i\Xvec_s$,
$\Xvec_i\Xvec_t$ and $\Yvec_k\Yvec_{\ell}$ again reduces the cost. A similar argument applies to $i \neq j$ and $k=\ell$. Finally ($iii$) if $i=j$ and $k=\ell$, there are two edges between $\Xvec_i$ and $\Yvec_k$. In this case, choose two edges $\Xvec_s\Xvec_t$ and 
$\Yvec_u\Yvec_v$, where $s\neq t\neq i$ and $u\neq v\neq k$. Now, replacing $\Xvec_i\Yvec_k$ (both copies), $\Xvec_s\Xvec_t$ and $\Yvec_u\Yvec_v$ by
$\Xvec_i\Xvec_s$,
$\Xvec_i\Xvec_t$,  $\Yvec_k\Yvec_u$ and $\Yvec_k\Yvec_v$ reduces the total derangement cost.

So, in each of the cases, the optimal derangement should have the minimum number of cross-sample edges. 
Consequently, $D_{m,n}$ and $D_{m,n}^{\circ}$ attain their respective minimum values 0 and 2, respectively, with probability
tending to 1 as $N$ and $d$ grow to infinity. Now, the consistency of the proposed tests can be proved using the fact that $
\mathbb{P}_{\mathrm{perm}}(D_{m,n}=0)<\alpha$
and $\mathbb{P}_{\mathrm{perm}}(D_{m,n}^{\circ}=2)<\alpha
$ (follows from the arguments used in the proof of Theorem 3.1).
\end{proof}

\begin{proof}[{\bf Proof of Theorem 3.2}] Write $a_d=S_{FF,d}$, $b_d=S_{GG,d}$ and $c_d=S_{FG,d}$ so that $S_d=2c_d-a_d-b_d$. Under the assumption $\liminf_{d\to\infty}d^{-1}\operatorname{tr}(\Sigma_i)>0$ for $i=1,2$, the limiting values of $a_d$, $b_d$ and $c_d$ turn out to be positive. So, there exists a constant $c_0>0$ such that for all sufficiently large $d$, $\min\{a_d,b_d,c_d\}\ge c_0>0$.  Set $\eta={c_0}/(16N)$ and let ${\cal A}_d$ be the event that simultaneously for all relevant pairs 
\begin{align*}
    \left| d^{-1}{\|\Xvec_i-\Xvec_j\|^2} -E\left\{d^{-1}{\|\Xvec_i-\Xvec_j\|^2}\right\} \right| & < \eta S_d, ~~~~ \left| d^{-1}{\|\Yvec_i-\Yvec_j\|^2} -E\left\{d^{-1}{\|\Yvec_i-\Yvec_j\|^2}\right\} \right| < \eta S_d\\ & \mbox{and   } \left| d^{-1}{\|\Xvec_i-\Yvec_j\|^2} -E\left\{d^{-1}{\|\Xvec_i-\Yvec_j\|^2}\right\} \right| < \eta S_d. 
\end{align*}  By Chebyshev's inequality and from the definition of $V_d$, for every relevant pair (${\bf U},{\bf V})$, we have \[ P\left( \left| d^{-1}{\|{\bf U}-{\bf V}\|^2} -E\left\{d^{-1}{\|{\bf U}-{\bf V}\|^2}\right\} \right|\ge \eta S_d \right) \le {V_d^2}/({\eta^2S_d^2}). \] Since $V_d/S_d \rightarrow 0$ as $d \rightarrow \infty$ and the sample sizes $m$, $n\geq 2$ are finite (i.e., the number of relevant pairs is also finite), by the union bound, we get $\P({\cal A}_d^c) \le \binom{N}{2}\frac{1}{\eta^2}\big({V_d}/{S_d}\big)^2 \longrightarrow 0$ as $d \rightarrow \infty$.

Using
$
|\sqrt{x}-\sqrt{y}|
=
\frac{|x-y|}{\sqrt{x}+\sqrt{y}}
$
and $\min\{a_d,b_d,c_d\}\ge c_0$, on ${\cal A}_d$ we obtain,
simultaneously for all relevant pairs,
$
\left|
{\|\Xvec_i-\Xvec_j\|}/{\sqrt d}-a_d
\right|
<
{S_d}/{16N},
$
$
\left|
{\|\Yvec_i-\Yvec_j\|}/{\sqrt d}-b_d
\right|
<
{S_d}/{16N},
$
and
$
\left|
{\|\Xvec_i-\Yvec_j\|}/{\sqrt d}-c_d
\right|
<
{S_d}/{16N}.
$
Let
$
\epsilon_d={S_d}/{16N}.
$
Now, consider a derangement $\pi$ which has $2t$ cross-sample edges. For such a derangement, on ${\cal A}_d$, $$\mathcal{C}_d(\pi)/\sqrt{d}\geq (m-t)(a_d-\epsilon_d)+(n-t)(b_d-\epsilon_d)+2t(c_d-\epsilon_d)=ma_d+nb_d+tS_d-N\epsilon_d,$$ 
where $\epsilon_d=S_d/(16N)$. Note that if we choose an arbitrary derangement $\pi_0$ with $0$ cross-sample edges,  on $\mathcal{A}_d$, we have ${\mathcal C}_d(\pi_0)/\sqrt{d}\leq m(a_d+\epsilon_d)+n(b_d+\epsilon_d) = ma_d+nb_d+N\epsilon_d$. 
Hence , $\big(\mathcal{C}_d(\pi)-\mathcal{C}_d(\pi_0)\big)/\sqrt{d}\geq tS_d-2N\epsilon_d
>0,$ So, for large $d$, for an optimal derangement, we must have $t=0$. This implies $\P(D_{m,n}=0)\rightarrow 1$ as $d\rightarrow\infty$. 

Note that any Hamiltonian cycle contains at least two cross-sample edges. Proceeding in the same way as above, we can show that for large $d$, on ${\cal A}_d$, an optimal derangement (in the class of derangements containing a Hamiltonian cycle) cannot have more than two cross sample edges. So, for ODER-SHC test, $\P(D_{m,n}^{\circ}=2) \rightarrow 1$ as $d \rightarrow \infty$.

Now, the proof can be completed by showing $\mathbb{P}_{\mathrm{perm}}(D_{m,n}=0) <\alpha$ and $\mathbb{P}_{\mathrm{perm}}(D_{m,n}^{\circ}=2) <\alpha$ under the given condition. This can be shown following the arguments used in the proof of Theorem 3.1. 
\end{proof}

\begin{rem}
Using the same notations as in the proof of Theorem~\ref{thm:joint-consistency}, we have $S_d=2c_d-a_d-b_d$, where
$
a_d=S_{FF,d}$, $
b_d=S_{GG,d}$ and $
c_d=S_{FG,d}$. Define $s_d=\min\{a_d,b_d\}$ and note that
$
c_d-s_d
=
({S_d+|a_d-b_d|})/{2}
\geq{S_d}/{2}$. Also  define  $M_{11}=\max_{1\le i<j \le m}\|\Xvec_i-\Xvec_j\|$, $M_{22}=\max_{1\le i<j \le n}\|\Yvec_i-\Yvec_j\|$ and $M_{12}=\min_{i,j}\|\Xvec_i-\Yvec_j\|$.
On the event ${\cal A}_d$ (as defined in the proof of Theorem 3.2), we have

$$\left(M_{12}-\min\{M_{11},M_{22}\}\right)/\sqrt{d}
\ge \left(c_d-\frac{S_d}{16N}\right)
-
\left(s_d+\frac{S_d}{16N}\right)
\geq
\frac{S_d}{2}-\frac{S_d}{8N}
>0.
$$
Under the condition of Theorem 3.2, since $\P({\cal A}_d) \rightarrow 1$ as $d \rightarrow \infty$, and on ${\cal A}_d$, 
$M_{12}>\min\{M_{11},M_{22}\}$,  either all $\Xvec\Xvec$-type edges or all $\Yvec\Yvec$-type edges turn out to be smaller than all $\Xvec\Yvec$-type edges with probability tending to $1$ as $d$ grows to infinity. So, following the argument of Remark 1, one can prove the high-dimensional consistency of the implemented version of the ODER-SHC test where the SHC is approximated using the cheapest-link algorithm. 
\end{rem}

\begin{lem}
\label{lem:sparse-rates-L2}
Suppose $\Xvec_1,\Xvec_2,\ldots,\Xvec_m$ are independent copies of $\Xvec \sim F$, and $\Yvec_1,\Yvec_2,\ldots,\Yvec_n$ are independent copies of $\Yvec \sim G$, where $F$ and $G$ satisfy (A1)-(A3).  Also assume that $\lim_{d \rightarrow \infty} d^{-1}{\text tr}(\sigmat_i)>0$ for $i=1,2$. If 
$\|\muvec_1-\muvec_2\|^2\gg d^{1/2}$
or
$\left|\operatorname{tr}(\sigmat_1)-\operatorname{tr}(\sigmat_2)\right|\gg d^{3/4}$, then $d^{1/2}S_d\rightarrow \infty $ as $d \rightarrow \infty$.
\end{lem}

\begin{proof}

Let $a_d=S_{FF,d}$,
$b_d=S_{GG,d}$,
and
$c_d=S_{FG,d}$ so that
$S_d=2c_d-a_d-b_d$. Also define $\nu_d^2=\|\muvec_1-\muvec_2\|^2/d$,  $\sigma_{1d}^2=\text{tr}(\sigmat_1)/d$ and $\sigma_{2d}^2=\text{tr}(\sigmat_2)/d$. 
Using
$a_d=\sqrt{2}\sigma_{1d}$,
$b_d=\sqrt{2}\sigma_{2d}$,
and
$c_d=\sqrt{\sigma_{1d}^2+\sigma_{2d}^2+\nu_d^2}$,
we have 
$S_d(2c_d+a_d+b_d)=2\{(\sigma_{1d}-\sigma_{2d})^2+2\nu_d^2\}$.
Since the $\sigma_{id}$'s are bounded away from $0$ for large $d$, under Assumption (A3), we have $\sigma_{1d}+\sigma_{2d}\asymp 1$ and $2c_d+a_d+b_d\asymp1$.

Now 
$\|\muvec_1-\muvec_2\|^2\gg d^{1/2}$
implies $\nu_d^2\gg d^{-1/2}$, and therefore $S_d\gg d^{-1/2}$
or $d^{1/2}S_d \rightarrow \infty$ as $d \rightarrow \infty$. 

Also,
$\left|\operatorname{tr}(\sigmat_1)-\operatorname{tr}(\sigmat_2)\right|\gg d^{3/4}$
implies $|\sigma_{1d}^2-\sigma_{2d}^2|\gg d^{-1/4}$.
Since $\sigma_{1d}+\sigma_{2d}\asymp 1$,
$|\sigma_{1d}-\sigma_{2d}|\asymp |\sigma^2_{1d}-\sigma^2_{2d}|\gg d^{-1/4}$ or
$(\sigma_{1d}-\sigma_{2d})^2\gg d^{-1/2}$. 
This implies $S_d \gg d^{-1/2}$ or $d^{1/2}S_d \rightarrow \infty$ as $d \rightarrow \infty$.
\end{proof}

\begin{proof}[{\bf Proof of Theorem~\ref{thm:HDLSS-L1}}]
Define
$A_d=\E\|\Xvec_1-\Xvec_2\|_1/d,$
$B_d=\E\|\Yvec_1-\Yvec_2\|_1/d$,
$C_d=\E\|\Xvec_1-\Yvec_1\|_1/d$ and $\Delta_d^{(1)}=2C_d-A_d-B_d$. Under Assumptions (A1$^\circ$) and (A2$^\circ$), as $d\rightarrow \infty$, we have
\begin{align*}
    \left|d^{-1}{\|\Xvec_i-\Xvec_j\|_1} -A_d\right|&
\stackrel{P}{\longrightarrow}0 \mbox{  for  }  1\le i<j \le m, 
\left|d^{-1}{\|\Yvec_i-\Yvec_j\|_1}-B_d\right|
\stackrel{P}{\longrightarrow}0 \mbox{  for  } 1\le i<j \le n\\
& \mbox{and   }
\Big|d^{-1}{\|\Xvec_i-\Yvec_j\|_1}
-C_d\Big|\stackrel{P}{\longrightarrow}0 \mbox{  for   } 1\le i\le m,~ 1\le j\le n.\end{align*}

For a derangement $\pi$, let
$T(\pi)$ be the number of directed cross-sample edges and 
$\mathcal C^{(1)}_d(\pi)
=
\sum_{i=1}^{m+n}\|Z_i-Z_{\pi(i)}\|_1$
be the total cost.
Since $m$ , $n\geq 2$ are fixed,
one can check that
\begin{align*}
\frac{\mathcal C_d^{(1)}(\pi)}{d}
~=~
\left(m-\frac{T(\pi)}2\right)A_d
+
\left(n-\frac{T(\pi)}2\right)B_d
+
T(\pi)C_d
+
o_P(1)
~ =~ mA_d + nB_d + T({\pi)} \Delta_d^{(1)}/2 + o_P(1).
\end{align*}

Since the number of admissible derangements is finite, the above \(o_P(1)\) term can be taken uniformly over \(\pi\).\\
 Note that
$\Delta_d^{(1)}
=
\frac1d\sum_{q=1}^d
\left[
2\E|X_1^{(q)}-Y_1^{(q)}|
-
\E|X_1^{(q)}-X_2^{(q)}|
-
\E|Y_1^{(q)}-Y_2^{(q)}|
\right]=\frac{1}{d}\sum_{q=1}^{d}{\cal E}_{F,G}^{(q)}.$ So, under the given condition $\liminf_{d \rightarrow \infty}\Delta_d^{(1)}>0$.
For the optimal derangement, we have
$\P(D_{m,n}=0)\longrightarrow1$
as $d$ grows to infinity. If we restrict ourselves to derangements having only one cycle of size $N$, the minimum number of cross-matches turns out to be $2$. So, in that case, we have 
$\P(D_{m,n}^{\circ}=2) \rightarrow 1$ as $d \rightarrow \infty$. Now, the rest can be proved by showing $\mathbb{P}_{\mathrm{perm}}(D_{m,n}=0)<\alpha$ and $\mathbb{P}_{\mathrm{perm}}(D_{m,n}^{\circ}=2)<\alpha$
as in the proof of Theorem 3.1.
\end{proof}

\begin{rem} 
Let $A$, $B$ and $C$ be the limiting values of $
A_d=d^{-1}E\|\Xvec_1-\Xvec_2\|_1$, 
$B_d=d^{-1}E\|\Yvec_1-\Yvec_2\|_1
$ and 
$
C_d=\frac{1}{d}\E\|\Xvec_1-\Yvec_1\|_1$, respectively. Under the condition of  Theorem~\ref{thm:HDLSS-L1}, we have
$
\lim_{d \rightarrow \infty}
d^{-1}\sum_{q=1}^d{\cal E}_{F,G}^{(q)} = 2C-A-B
>0,
$ which implies $C>\min\{A,B\}$. So, the high-dimensional consistency of the implemented version of the ODER-SHC-L1 test based on the cheapest-link algorithm follows using the same arguments as in Remark 1.


\end{rem}


\begin{proof}[{\bf Proof of Theorem~\ref{thm:HDLSS-L1-sparse}}]

Let
$A_d=d^{-1}\E\|\Xvec_1-\Xvec_2\|_1$,
$B_d=d^{-1}\E\|\Yvec_1-\Yvec_2\|_1$
and
$C_d=d^{-1}\E\|\Xvec_1-\Yvec_1\|_1$, so that
$S_d^{\circ}=2C_d-A_d-B_d$.
Set $\eta^{\circ}_d=S_d^{\circ}/16N$ and let ${\mathcal A}^{\circ}_d$ denote the event on which, simultaneously for all relevant pairs,
$\left|d^{-1}\|\Xvec_i-\Xvec_j\|_1-A_d\right|< \eta^{\circ}_d$,
$\left|d^{-1}\|\Yvec_i-\Yvec_j\|_1-B_d\right|< \eta^{\circ}_d$
and
$\left|d^{-1}\|\Xvec_i-\Yvec_j\|_1-C_d\right|< \eta^{\circ}_d$.

For each pair (${\bf U}, {\bf V})$, using Chebyshev's inequality, one gets
$\P\left(\left|d^{-1}\|U-V\|_1-\E\{d^{-1}\|U-V\|_1\}\right|\ge \eta^{\circ}_d\right)\leq(V_d^{\circ}/\eta^{\circ}_d)^2$.
So, a union bound over all relevant pairs yields
$\P(({\mathcal A}_d^{\circ})^ c)\leq \binom{N}{2} (16N)^2(V_d^{\circ}/S_d^{\circ})^2$, and hence $\P(({\mathcal A}_d^{\circ})^ c)
\rightarrow 0$ as $d \rightarrow \infty$
(since $N$ is finite and $V_d^{\circ}/S_d^{\circ}\rightarrow 0$ as $d \rightarrow \infty$).


Now, following the arguments as in the proof of Theorem 3.2, one can show that $\P(D_{m,n}=0)\rightarrow 1$ and  $\P(D_{m,n}^{\circ}=2)\rightarrow 1$ as $d \rightarrow \infty$.
The rest of the proof follows from the arguments given in the proof of Theorem~\ref{thm:HDLSS-derangement} regarding the calculation of these two probabilities under the permutation distribution. 
\end{proof}

\begin{rem} 
Following the notations used in the proof of Theorem~\ref{thm:HDLSS-L1-sparse}, we have
$S_d^{\circ}=2C_d-A_d-B_d$, where
$A_d={d}^{-1}\E\|\Xvec_1-\Xvec_2\|_1$, 
$B_d=d^{-1}\E\|\Yvec_1-\Yvec_2\|_1,
$
and $
C_d=d^{-1}\E\|\Xvec_1-\Yvec_1\|_1$.
Also define $s_d^{\circ}=\min\{A_d,B_d\}$. Note that on ${\cal A}_d^{\circ}$ (as defined in the proof of Theorem~\ref{thm:HDLSS-L1-sparse}), 
$
C_d-s_d^{\circ}
=
\big({S_d^{\circ}+|A_d-B_d|}\big)/{2}
\geq {S_d^{\circ}}/{2}>0.
$ 
Since $\P({\cal A}_d^{\circ}) \rightarrow 1$ as $d \rightarrow \infty$, under the conditions of Theorem 3.4, the high-dimensional consistency of the implemented version of the ODER-SHC-L1 test follows using the same argument as in Remark 2.
\end{rem}

\subsection{Results on the large-sample behavior of the derangement-based tests}
\vspace{0.1in}

To prove the asymptotic normality of $D_{m,n}$ stated in Theorem 6.1, we use the following two lemmas. The first one is taken from \cite{ross2011fundamentals}, and the proof is available in that article.

\begin{lem} {\bf (Theorem 3.7 in \cite{ross2011fundamentals}).}
\label{lem:exchangeable-pair}
Let $(W,W')$ be an exchangeable pair such that $\E W=0$, $\E W^2=1$, and
$\E(W'-W\mid W)=-\lambda W$  for some $\lambda\in(0,1)$.
If $Z\sim N(0,1)$, then
$$d_W(W,Z)=\sup_{\|h'\|_\infty\leq1}|\E h(W)-\E h(Z)|~\leq\frac{1}{\sqrt{2\pi\lambda}}\sqrt{\Var\{\E[(W'-W)^2\mid W]\}}+\frac{1}{3\lambda}\E|W'-W|^3.$$
\end{lem}


\begin{lem}
\label{lem:slice-poincare}
Let
$\mathcal S_{N,m}=\{\ell=(\ell_1,\ldots,\ell_N)\in\{0,1\}^N:\sum_{i=1}^N\ell_i=m\}$,
and let $L$ be uniformly distributed on $\mathcal S_{N,m}$. Choose a random ordered pair $(I,J)$ from a uniform distribution on $S_0=\{(i,j): 1\le i,j \le N, i\neq j\}$, independently of $L$, to get $L'$ by interchanging $I$-th and $J$-th elements of $L$. Then, for every $H:\mathcal S_{N,m}\rightarrow\mathbb R$,
$$\Var\{H(L)\}\leq\frac{N-1}{4}\E\{H(L')-H(L)\}^2.$$
\end{lem}

\begin{proof}
Since $L$ is uniformly distributed on $\mathcal S_{N,m}$, the random-transposition Poincaré inequality \citep[see Lemma~1.3 in][]{salez2021sharp} gives
$
\operatorname{Var}{H(L)}
\leq
\frac{1}{2N}
\sum_{1\leq i<j\leq N}
E\left[\{H(L_{ij})-H(L)\}^{2}\right],
$
where $L_{ij}$ is obtained from $L$ by interchanging its $i$-th and $j$-th coordinates. Since $(I,J)$ is uniformly distributed over the $N(N-1)$ ordered pairs in $S_0$ and $L_{ij}=L_{ji}$,
\[
E\left[\{H(L')-H(L)\}^{2}\right]
=
\frac{1}{N(N-1)}
\sum_{i\ne j}
E\left[\{H(L_{ij})-H(L)\}^{2}\right]
=
\frac{2}{N(N-1)}
\sum_{1\le i<j\le N}
E\left[\{H(L_{ij})-H(L)\}^{2}\right].
\]
Combining the above inequalities, one gets
$
\operatorname{Var}\{H(L)\}
\leq
\frac{N-1}{4}
E\left[\{H(L')-H(L)\}^{2}\right].
$
\end{proof}

\begin{proof}[{\bf Proof of Theorem~\ref{thm:null-asymptotic}}]

Given the pooled sample $\Zvec_1,\Zvec_2,\ldots,\Zvec_N$ and the selected derangement $\pi_N$, the only randomness comes from permutation of $m$ $\Xvec$-labels and $n=N-m$ $\Yvec$-labels. The conditional mean $\nu_{m,n}$ and variance $\omega_{m,n}^2$ of the permutation distribution of $D_{m,n}$ have already been derived in \eqref{eq:null-mean} and \eqref{eq:null-variance-simplified}.
Taking $p_N=m/N$, $q_N=n/N$, and $\kappa_N=C_{N}/N$, from \eqref{eq:null-variance-simplified}, one gets
$\omega_{m,n}^2/N=4p_N^2q_N^2(1+2\kappa_N)+o(1)$. Since $p_N\rightarrow p\in(0,1)$ and $\kappa_N$ is bounded between $0$ and $1$, there exist constants $c_1,c_2>0$ such that
$c_1N\leq\omega_{m,n}^2\leq c_2N$
for all sufficiently large $N$.

Let us define $w_{ij}$ as the number of edges between $\Zvec_i$ and $\Zvec_j$ in the optimal derangement graph. Clearly, $w_{ij}=w_{ji}$ and $0\leq w_{ij}\leq2$ for all $1\le i<j \le N$.
 Also, corresponding to a permutation,
define $L$ as an $N$-dimensional
vector, whose $i$-th element ($i=1,2,\ldots,N$) $L_i$ is $1$ if
an $\Xvec$-label is assigned to $\Zvec_i$, $0$ otherwise. 

Now, note that the test statistic $D_{m,n}$ can be expressed as $D_{m,n}=\sum_{1\leq i<j\leq N}w_{ij}I\{L_i\neq L_j\}$. Choose an ordered pair $(I,J)$ randomly from the uniform distribution on $S_0=\{(i,j): 1\le i,j \le N, i\neq j\}$, and interchange the corresponding elements of $L$ to get a random vector $L'$ and the resulting statistic $D_{m,n}'$. 
Since the permutation distribution is invariant under such interchanging, $(D_{m,n}, D_{m,n}')$ is exchangeable. For any fixed $i$ and $j$ ($1\le i,j\le N$), 
let $L^{(ij)}$ be the vector obtained from $L$ by interchanging $i$-th and $j$-th elements, and $D_{m,n}^{(ij)}$ be the corresponding statistic. Clearly, 
$D_{m,n}^{(ij)}=D_{m,n}$ if $L_i=L_j$. Otherwise, one can show that
$$D_{m,n}^{(ij)}-D_{m,n}=4-2(c_i+c_j)+2w_{ij}, \mbox{  where   }
c_i={\textstyle \sum_{k=1}^Nw_{ik}I\{L_i\neq L_k\}}.$$
Let $A=\{i:L_i=1\}$ and $B=\{i:L_i=0\}$. Since each cross-sample edge contributes its weight once to each of its two endpoint classes,
$\sum_{i\in A}c_i=\sum_{j\in B}c_j=D_{m,n}$
and
$\sum_{i\in A}\sum_{j\in B}w_{ij}=D_{m,n}$.
Therefore
$\sum_{i\in A}\sum_{j\in B}(D_{m,n}^{(ij)}-D_{m,n})=4mn-2(N-1)D_{m,n}$.
The same expression is obtained over $B\times A$, and hence
$$\E(D_{m,n}'-D_{m,n}\mid L_1,\ldots,L_N)=-\frac4N(D_{m,n}-\nu_{m,n}).$$
Now, take
$W_N=(D_{m,n}-\nu_{m,n})/\omega_{m,n}$
and
$W_N'=(D_{m,n}'-\nu_{m,n})/\omega_{m,n}$. Clearly $(W_N,W_N')$ forms an exchangeable pair with
$\E(W_N'-W_N\mid W_N)=-\lambda_NW_N$
with
$\lambda_N=4/N$. Let us now look at the two terms on the right side of the inequality stated in Lemma~\ref{lem:exchangeable-pair}. Since $\sum_jw_{ij}=2$ and $0\leq w_{ij}\leq2$,
$|D_{m,n}'-D_{m,n}|\leq4$.
Consequently,
$$\frac{1}{3\lambda_N}\E|W_N'-W_N|^3\leq\frac{64}{3\lambda_N\omega_N^3}=O(N^{-1/2}).$$

For $\ell\in\mathcal S_{N,m}$, define
$G_N(\ell)=\E[(D_{m,n}'-D_{m,n})^2\mid L=\ell]$.
For fixed $i\neq j$,  $D_{m,n}^{(ij)}-D_{m,n}$  depends only on the labels of $\Zvec_i$, $\Zvec_j$ and their neighboring (adjacent) vertices. Since every vertex has at most two distinct neighbors, changing the labels of two observations can affect at most $O(N)$ of the $N(N-1)$ terms in the expression of $G_N$, and in each case, the squared increment is bounded by $16$. Thus there exists a constant $C$ such that
$|G_N(\ell')-G_N(\ell)|\leq C/N$
whenever $\ell'$ is obtained from $\ell$ by interchanging the labels at two locations.

Therefore, Lemma~\ref{lem:slice-poincare} gives
$\Var\{G_N(L)\}=O(N^{-1})$.
Since
$\E[(W_N'-W_N)^2\mid L]=G_N(L)/\omega_{m,n}^2$,
we obtain
$\Var\{\E[(W_N'-W_N)^2\mid W_N]\}\leq\Var\{G_N(L)/\omega_{m,n}^2\}=O(N^{-3})$. So,
it follows that
$$\frac{1}{\sqrt{2\pi}\lambda_N}\sqrt{\Var\{\E[(W_N'-W_N)^2\mid W_N]\}}=O(N^{-1/2}).$$

Now, from Lemma~\ref{lem:exchangeable-pair}, one gets
$d_W(W_N,Z)\rightarrow0$
for $Z\sim N(0,1)$, and hence, under the given conditions
$W_N=({D_{m,n}-\nu_{m,n}})/{\omega_{m,n}}\xrightarrow{d}N(0,1)$ as $N$ tends to infinity.
\end{proof}

\begin{rem}
Note that the arguments given in the proof of Theorem 6.1 depend only on the structure of the selected derangement graph, and hence they hold both for $D_{m,n}$ and $D_{m,n}^{\circ}$. The only difference is that the expressions of $\omega_{m,n}$ in these two cases are slightly different, though they are of the same asymptotic order. Moreover, the proof remains the same irrespective of the derangement cost function. So, asymptotic normality of the null distributions of $D_{m,n}$ and $D_{m,n}^\circ$ holds even when the $\ell_1$ distance is used as the derangement cost function.
\end{rem}

\begin{lem}
\label{lem:hamiltonian-path}
For any $N\ge 2$ points in $[-R,R]^d$, there exists a Hamiltonian path ${\mathbb H}_N$ such that 
${\cal L}({\mathbb H}_N) \le C_d R N^{1-1/d}$, where ${\cal L}({\mathbb H}_N)$ denotes the total length of ${\mathbb H}_N$ and $C_d$
is a positive constant that depends only on $d$.
\end{lem}

\begin{proof}
    Let $\Zvec_1,\ldots,\Zvec_N$ be $N$ ($N\ge 2$) points in $[-R,R]^d$. Take
$M=\lceil N^{1/d}\rceil$ and
partition the hypercube $[-R,R]^d$ into $M^d$ congruent sub-cubes, each having side length
$h=2R/M$ so that the diameter (in terms of Euclidean distance) is at most $\sqrt d h$.
Order the $M^d$ sub-cubes in a serpentine fashion so that two consecutive sub-cubes share a face. Let
$Q_{i_1},\ldots,Q_{i_s}$ be the non-empty sub-cubes arranged in that order and let $n_\ell$ be the number of points in
$Q_{i_\ell}$. Within each occupied sub-cube $Q_{i_\ell}$, join its points in an arbitrary order to get a path ${\cal P}_{i_\ell}$.
Since every two points in the same sub-cube are at Euclidean distance at most
$\sqrt d\,h$, the total length of these $s$ paths is bounded by
$\sqrt d\,h\sum_{\ell=1}^s(n_\ell-1)\leq\sqrt d\,hN$.

Next, we connect these paths obtained in successive non-empty sub-cubes.
Since the centers of two consecutive sub-cubes in the serpentine order are at a distance $h$, the distance between the centres of $Q_{i_\ell}$ and $Q_{i_{\ell+1}}$ is at most $(i_{\ell+1}-i_{\ell})h$. Again, in any sub-cube, the distance between any point and the centre is at most $\sqrt{d} h/2$. So, using the triangle inequality, the distance from the last point of
${\cal P}_{i_\ell}$ to the first point of ${\cal P}_{i_{\ell+1}}$ is bounded by
$\sqrt{d}h\{i_{\ell+1}-i_\ell+1\}$.
Summing over $\ell$ and using
$\sum_{\ell=1}^{s-1}(i_{\ell+1}-i_\ell)=i_s-i_1\leq M^d$,
the total length of the Hamiltonian path is bounded by
$\sqrt{d}h(M^d+N)=\sqrt{d}\frac{2R}{M}(M^d+N) \le C_dRN^{1-1/d}$, where 
$C_d$ is a positive constant that depends only on $d$.
\end{proof}
\begin{lem}
\label{lem:hamiltonian-bound}
Let $d\geq2$ and let $\phi,\phi_1,\phi_2,\ldots$ be a sequence of densities with the same support such that
$\sup_{\{\phi>0\}}|(\phi_N/\phi
)-1|\rightarrow0$ as $N \rightarrow \infty$.
Suppose that for some $r>d/(d-1)$,
$\int_{\mathbb R^d}\|\zvec\|_2^r~\phi(\zvec)\,d\zvec<\infty$.
If
$\Zvec_{1,N},\ldots,\Zvec_{N,N}\stackrel{{i.i.d.}}{\sim}\phi_N$
and ${\cal L}_N^{\circ}$ is the total length of the SHC (in terms of Euclidean distance) through these points, then
$\E ({\cal L}_{N}^{\circ})=O(N^{1-1/d})$.
\end{lem}

\begin{proof}
Define
$S_0=\{z:\|z\|_\infty\leq1\}$
and 
$S_j=\{z:2^{j-1}<\|z\|_\infty\leq2^j\}$ for $j\ge 1$. 
Also define
$n_{j,N}=\sum_{i=1}^NI\{Z_{i,N}\in S_j\}$
and
$p_{j,N}=P_{\phi_N}(Z_{1,N}\in S_j)$. For every non-empty shell $S_j$, its points lie in $[-2^j,2^j]^d$. Using Lemma~\ref{lem:hamiltonian-path}, there exists a Hamiltonian path ${\mathbb H}_j$ on those points with length at most
$C_d2^jn_{j,N}^{\gamma}$ for $\gamma=1-1/d$. 
Let $S_{i_1},S_{i_2},\ldots,S_{i_k}$ be the non-empty shells arranged in serpentine order. Connecting one end-point of ${\mathbb H}_{i_\ell}$ with an end-point of ${\mathbb H}_{i_{\ell+1}}$ for $\ell=1,2,\ldots,k-1$, we get a Hamiltonian path on $N$ vertices.
Note that the distance between a point on ${\mathbb H}_{i_{\ell}}$ and a point on  ${\mathbb H}_{i_{\ell+1}}$ is at most $C_d2^{i_{\ell+1}}$. Finally, the Hamiltonian cycle is formed by joining the open ends of ${\mathbb H}_{i_1}$ and ${\mathbb H}_{i_k}$ by an edge of length bounded by $C_d 2^{i_k}$. Therefore, ${\cal L}_N^{\circ} \le C_d\sum_{j\geq0}2^j\{n_{j,N}^{\gamma}+I(n_{j,N}>0)\}$.

Since $x^\gamma$ is concave and
$\E n_{j,N}=Np_{j,N}$, we have 
$\E n_{j,N}^{\gamma}\leq(Np_{j,N})^\gamma$.
Moreover, for $j\geq1$, on the event $\{\Zvec_{1,N}\in S_j\}$ we have
$\|\Zvec_{1,N}\|_2\geq\|\Zvec_{1,N}\|_\infty>2^{j-1}$.
The uniform $r$th-moment bound and Markov's inequality therefore give
$p_{j,N}\leq \E{\|\Zvec_1\|^r}/2^{(j-1)r} \le C_02^{-rj}$ for some $C_0>0$. Consequently,
$$\E\sum_{j\geq0}2^jn_{j,N}^{\gamma}
\leq C_0^{\gamma}N^\gamma
\Big\{1+\sum_{j\geq1}2^{j-rj\gamma }\Big\}
=C_0N^\gamma
\Big\{1+\sum_{j\geq1}2^{-j(r\gamma-1)}\Big\} ~\mbox{  for some   } C_0>0.$$
Since $r>d/(d-1)=1/\gamma$, we have $r\gamma>1$, and hence
$\E\sum_{j\geq0}2^jn_{j,N}^{\gamma}=O(N^\gamma)$.

Also,
$\P(n_{j,N}>0)=1-(1-p_{j,N})^N \leq\min\{1,Np_{j,N}\}$.
Choose $J_N$ such that
$2^{rJ_N}\leq N<2^{r(J_N+1)}$.
Then
\begin{align*}
\E\sum_{j\geq0}2^jI(n_{j,N}>0) =\sum_{j\geq0} 2^{j}\{1-(1-p_{j,N})^N\} \le  \sum_{j\leq J_N}2^j+ \sum_{j>J_N}2^{j} Np_{j,N}
\end{align*}
Now, $\sum_{j\leq J_N}2^j=O(2^{J_N})=O(N^{1/r})$,
while
$N\sum_{j>J_N}2^jp_{j,N}
\leq
C_0N\sum_{j>J_N}2^{-(r-1)j}
=\frac{C_0N}{2^{rj_N+1}}2^{j_N+1}\leq c_02^{j_N+1}=O(N^{1/r})$. 
Since $1/r<\gamma$,
$N^{1/r}=o(N^\gamma)$.

Combining the two bounds, we get
$\E {\cal L}_{N}^{\circ}
=O(N^\gamma)
=O(N^{1-1/d})$.
\end{proof}

\begin{proof}[{\bf Proof of Theorem~\ref{thm:large-sample-consistency}}]

Define
$h(z)=pf(z)+qg(z)$
and
$h_N(z)=p_Nf(z)+q_Ng(z)$, where $p+q=1$, $p_N=m/N$, and $q_N=1-p_N$.
On $\{h>0\}$,
$$\left|\frac{h_N(z)}{h(z)}-1\right|=|p_N-p|~\frac{|f(z)-g(z)|}{pf(z)+qg(z)}\leq|p_N-p|~\left(\frac{1}{p}+\frac{1}{q}\right).$$
Since $|p_N-p| \rightarrow 0$, we have $\sup_{\{h>0\}}|h_N/h-1|\rightarrow0$ as $N \rightarrow \infty$.

Let ${\cal L}_N^{\circ}$ be the length of the SHC (i.e., optimal derangement graph with one cycle of size $N$) when the Euclidean distance is used as the transportation cost. From Lemma A.6, we know that $\E({\cal L}_N^{\circ})=O(N^\gamma)$, where $\gamma=(d-1)/d$.

Note that in the SHC, every vertex has exactly one outgoing edge and one incoming edge. 
Now, for $a>0$, define
$M_{a,N}=\sum_{i=1}^NI\{\|Z_{i,N}-Z_{\pi_N(i),N}\|>aN^{-1/d}\}$.
Since every non-zero term in $M_{a,N}$ contributes more than $aN^{-1/d}$ to ${\cal L}_N^{\circ}$, 
we have 
$aN^{-1/d}~\E M_{a,N}\leq \E({\cal L}_N^{\circ})$. This implies $\E M_{a,N} \le CN/a$, for some positive constant $C$.
So, under permutation equivariance, we have
$\P\{\|Z_{1,N}-Z_{\widehat\pi_N(1),N}\|_2>aN^{-1/d}\}=N^{-1}\E M_{a,N}\leq C/a$.
Therefore
$\lim_{a\to\infty}\limsup_{N\to\infty}P\{\|Z_{1,N}-Z_{\widehat\pi_N(1),N}\|_2>aN^{-1/d}\}=0$.

So, the SHC satisfies the short-edge property mentioned in \cite{Arias-Castro2016}.
Moreover, every vertex of a Hamiltonian cycle has one incoming and one outgoing edge. So, it also satisfies the bounded-degree condition. Hence, from Theorem~1 of ~\cite{Arias-Castro2016},
%
$${D_{m,n}^{\circ}}/{N}\stackrel{a.s.}{\longrightarrow} 2pq\int_{\mathbb R^d}\frac{f(z)g(z)}{pf(z)+qg(z)}\,dz:=\Lambda_p(F,G).$$
So, under ${\cal H}_0:F=G$, $D_{m,n}^{\circ}/N$ converges to $2pq$ almost surely. 
One can show that
$$2pq-\Lambda_p(F,G)=2p^2q^2\int_{\mathbb R^d}\frac{(f(z)-g(z))^2}{pf(z)+qg(z)}\,dz\ge 0,$$ where the equality holds if and only if $f=g$ almost everywhere.
So, the limiting value under ${\cal H}_1:F\neq G$ is lower than that under ${\cal H}_0$. Since the ODER-SHC test rejects ${\cal H}_0$ for small values of $D_{m,n}^{\circ}$ or $D_{m,n}^{\circ}/N$, under ${\cal H}_1$, this rejection probability converges to $1$ as $N$ tends to infinity. This proves the large-sample consistency of the ODER-SHC test.

The large-sample consistency of the ODER test also follows from a similar set of arguments. Let $\pi_N$ be the optimal derangement and ${\cal L}_N$ be the length of the corresponding optimal derangement graph. Clearly $\E({\cal L}_N)\le \E({\cal L}_N^{\circ})=O(N^{\gamma})$. So, following the same argument as above, we can show that the optimal derangement graph satisfies the short-edge property. Like the SHC, each vertex of the optimal derangement graph also has one incoming and one outgoing edge. Therefore,  from Theorem 1 of \cite{Arias-Castro2016}, $D_{m,n}/N \stackrel{a.s.}{\longrightarrow} \Lambda_p(F,G)$ as $N \rightarrow \infty$. 

Now, let us consider the permutation distribution of $D_{m,n}$. In the proof of Theorem ~\ref{thm:null-asymptotic}, the arguments are purely combinatorial and therefore remain valid under the alternative. This implies
that the permutation distribution of $D_{m,n}$ becomes degenerate at $2pq$ as $N$ tends to infinity.
Since $\Lambda_p(F,G)<2pq$ under any fixed alternative, this proves the large sample consistency of the ODER test. 
\end{proof}

\begin{rem}
    Since all norm-induced distances are equivalent in finite dimensions, for any norm-induced distance $\rho$, one can show that there exist $C_2\ge C_1> 0$ such that $C_1\|\zvec-\wvec\|\le \rho(\zvec,\wvec)\le C_2\|\zvec-\wvec\|$ for all $\zvec,\wvec$. Using this, under the condition
$\E_f(\rho(Z,0)^{r}) <\infty$ and $\E_g(\rho(Z,0)^{r}) <\infty$ for some $r>{d}/({d-1})$ (this is similar to Assumption (A5)), one can prove the short-edge property of 
the SHC and the optimal derangement graph based on $\rho$ as well. So, the large-sample consistency of the ODER and ODER-SHC test can be established using the same arguments as in the proof of Theorem 
\ref{thm:large-sample-consistency}.
\end{rem}

\begin{lem}
    Under the condition of Theorem~\ref{thm:large-sample-consistency}, the large-sample consistency holds for the ODER-SHC test even when the Hamiltonian cycle is constructed using the cheapest-link algorithm.
\end{lem}

\begin{proof}
Let ${\mathbb H}_N^{K}$ be the Hamiltonian cycle on $\Zvec_1,\Zvec_2,\ldots,\Zvec_N$ constructed using the cheapest-link algorithm. Here we shall show that ${\mathbb H}_N^{K}$ satisfies the short-edge property. The rest of the proof then follows using the same set of arguments as in the proof of Theorem 6.2. We shall prove this result for the test based on Euclidean distance only. But because of their equivalence in finite dimensions, this result holds for the tests based on any other norm-induced distance function.



First fix $a>0$ and define
$
s_N=a N^{-1/d}.
$ Consider the graph ${\cal G}_N(a)$ obtained by joining all relevant edges up to length $s_N$ following the cheapest-link algorithm.
Let
$
K_N(a)
$
be the number of connected components ${\cal G}_1,{\cal G}_2,\ldots,{\cal G}_{K_N(a)}$ in ${\cal G}_N(a)$, and
$
V_{\ell,N}$ be one of the two end points of ${\cal G}_{\ell}$ ($\ell=1,2,\ldots,K_N(a))$. Clearly, 
$\|V_{\ell,N}-V_{k,N}\|
>{s_N}$.

Now for any $\epsilon >0$, choose  $R>0$ such that $\P_{\phi}(\Zvec \notin [-R,R]^d)<\epsilon/2.$ 
Partition the hypercube $B_R=[-R,R]^d$ into sub-cubes having side length
$
\delta_N={s_N}/{\sqrt{d}}.
$
So, the diameter of each such sub-cube is
$
\sqrt d\,\delta_N
=
{s_N}.
$
Hence, no such
sub-cube can contain two distinct
$V_{j,N}$'s. 

The number of cubes required to cover $B_R$ is at most
$
\left(
1+\frac{2R}{\delta_N}
\right)^d
=
\left(
1+\frac{2\sqrt{d}R}{s_N}
\right)^d.
$
So,
$
K_N(a)
\leq
\left(
1+\frac{2\sqrt{d}R}{s_N}
\right)^d
+
\sum_{i=1}^N
I\{\|Z_{i,N}\|_\infty>R\}.
$
Now, taking expectations, dividing by $N$, and using
$s_N=a N^{-1/d}$, we obtain

\vspace{-0.1in}
$$
\frac{E K_N(a)}{N}
\leq
\frac{1}{N}
\left(
1+\frac{2\sqrt{d}R N^{1/d}}{a}
\right)^d
+
\Pr\{\|Z_{1,N}\|_\infty>R\} \le \left(
N^{-1/d}+\frac{2\sqrt{d}R}{a}
\right)^d
+ \epsilon/2.
\vspace{-0.1in}
$$
Hence, for every fixed $a>0$,
$
\limsup_{N\to\infty}
{E K_N(a)}/{N}
\leq
\left({2\sqrt{d}R}/{a}\right)^d
+
\epsilon/2.
$
So, for sufficiently large $a$, we have $\limsup_{N\to\infty}
{E K_N(a)}/{N} <\epsilon$. Since $\epsilon$ is arbitrary, this implies
$
\lim_{a\to\infty}
\limsup_{N\to\infty}
{E K_N(a)}/{N}
=0.
$

Note that the graph ${\cal G}_N(a)$ contains $N-K_N(a)$ edges, each of length at most $s_N$. Now, to form the Hamiltonian cycle, we need to add $K_N(a)$ edges. Therefore,
%
the number of edges of length bigger than $s_N$ in the Hamiltonian cycle is  at most $K_N(a)$, 
Hence, it follows that
$\lim_{a\to\infty}
\limsup_{N\to\infty}
\Pr\big\{
\|
Z_{1,N}
-
Z_{{\pi}^{\,K}_N(1),N}
\|_2
>
a N^{-1/d}
\big\}
\le 
\lim_{a\to\infty}
\limsup_{N\to\infty}
E K_N(a)/N =0,
$
where $\pi_N^K$ is the derangement map corresponding to the Hamiltonian cycle formed by the cheapest-link algorithm.
This proves the large-sample consistency of the resulting test.
\end{proof}

\bibliographystyle{apalike}
\bibliography{refs}

\end{document}